\documentclass[draft]{article}
\def\today{5.9.12} 
\usepackage{amsmath,amsfonts,amsthm,amssymb,amscd}
\renewcommand{\Re}{\mathop{\rm Re}\nolimits}
\renewcommand{\Im}{\mathop{\rm Im}\nolimits}

\theoremstyle{plain} \newtheorem{theorem}{Theorem}[section]
\newtheorem{lemma}[theorem]{Lemma}

\newtheorem{corollary}[theorem]{Corollary} \theoremstyle{definition}
\newtheorem{definition}[theorem]{Definition} \theoremstyle{remark}
\newtheorem{remark}[theorem]{Remark}

\newcommand{\R}{{\mathbb R}} \newcommand{\U}{{\mathcal U}}
\newcommand{\Hc}{{\mathcal H}}

\newcommand{\Z}{{\mathbb Z}}

\newcommand{\toro}{{\mathbb T}}

\newcommand{\e}{{\rm e}}

\newcommand{\Tr}{{\mathcal T}}

\newcommand{\Vc}{{\mathcal V}}

\newcommand{\Bc}{{\mathcal B}}

\newcommand{\I}{{\mathcal I}}
\newcommand{\G}{{\mathcal G}}

\newcommand{\resto}{{\mathcal R}}
\def\im{{\rm i}}

\newcommand{\Sc}{{\mathcal S}}
\newcommand{\Zc}{{\mathcal Z}}

\newcommand{\V}{{\mathcal V}}
\newcommand{\A}{{\mathcal A}}

\newcommand{\C}{\mathbb{C}}

\def\hs{{\frak h}}

\def\uno{{\kern+.3em {\rm 1} \kern -.22em {\rm l}}}

\def\der#1#2{\frac{\partial{#1}}{\partial {#2}}}
\def\di{{\rm d}}

\def\norma#1{\left\| #1\right\|}

\def\normas#1#2{N_{#1}^{\Sc}\left(#2\right)}
\def\normad#1#2{N_{#1}^{\nabla}\left(#2\right)}

\def\poisson#1#2{\left\{#1;#2\right\}}

\def\sleq{\preceq}

\def\nr{\normas{\tilde\delta_{r}}{\resto_{r}}}
\def\nz{\normad{\tilde\delta_{r}}{\Zc_{r}}}
\def\cg{{\cal C}}
\def\bl{{\bf l}}
\def\d{{\rm d}}
\def\e{{\rm e}}
\def\normaste#1#2#3{\norma{#1}_{L^{#2}_{\epsilon t}\bl^{#3}}}
\def\normast#1#2#3{\norma{#1}_{L^{#2}_{t}\bl^{#3}}}
\def\snorma#1{\left\langle\left|#1\right|\right\rangle}
\newcommand{\Ph}{\bl}

\numberwithin{equation}{section}

\begin{document}

\title{Asymptotic stability of breathers
in some Hamiltonian 
networks of weakly coupled oscillators}
\author{Dario Bambusi}

\date{\today}
\maketitle

\begin{abstract}
We consider a Hamiltonian chain of weakly coupled anharmonic
oscillators. It is well known that if the coupling is weak enough then
the system admits families of periodic solutions exponentially
localized in space (breathers). In this paper we prove asymptotic
stability in energy space of such solutions. The proof is based on two
steps: first we use canonical perturbation theory to put the system in
a suitable normal form in a neighborhood of the breather, second we
use dispersion in order to prove asymptotic stability. The main
limitation of the result rests in the fact that the nonlinear part of
the on site potential is required to have a zero of order 8 at the
origin. From a technical point of view the theory differs from that
developed for Hamiltonian PDEs due to the fact that the breather is
not a relative equilibrium of the system.
\end{abstract}

\section{Introduction}
\label{1} 

In this paper we consider the dynamical system with Hamiltonian
\begin{equation}
\label{0.1}
H:=\sum_{k\in\Z} \left[\frac{p_k^2+q_k^2}2+V(q_k)
  \right]+\frac\epsilon2 \sum_{k\in\Z}(q_{k+1}-q_{k})^2\ ,
\end{equation}
where $V$ is an analytic function having a zero of order at least 8 at
the origin. In 1994 MacKay and Aubry \cite{MA94} proved that if
$\epsilon$ is small enough, then there exist periodic solutions which
are exponentially localized in space (breathers). 

The problem of stability of the breathers has attracted a lot of
interest since the discovery of such objects and indeed linear stability has
been rapidly obtained through signature theory (see
\cite{MKS98}). Concerning the nonlinear stability, the only known
result for Hamiltonian networks ensures stability over times
exponentially long with $1/\epsilon$ \cite{Bam96}. However the
presence of dispersion suggests that the breathers should be
asymptotic stable (see e.g. \cite{Bam98}). (For nice
reviews on breathers see \cite{Aub97a,FlaWi98}.)

In the present paper we actually prove that breathers are
asymptotically stable, at least if the nonlinear part of the on site
potential fulfills $V(q)=O(|q|^8)$ as $q\to0$. More precisely we prove
that if the initial datum is close in the energy norm to a breather,
then the distance of the solution from the breathers, as a function of
time, is small as an element of $L^q_t(\R,\ell^r)$. As usual $(q,r)$
are admissible pairs (see eq.\eqref{a.p} below for a precise
definition).

We emphasize that such a result is one of the few examples of
asymptotic stability in Hamiltonian systems for object which are
neither equilibria nor relative equilibria. As far as we know the only
other known example is that of the solitary wave of the FPU system (see
\cite{FPII,FPIV,HofWay,Miz09,Miz11}). For the theory of asymptotic stability
of equilibria or relative equilibria see
e.g. \cite{SW90,BusPer92,Sig93,SofW99,Cuc01,GSNT04,BC09,Bam11a}.
\vskip10pt The proof consists of essentially 2 steps, the first one
consists in using canonical perturbation theory in order to put the
system in a suitable normal form. The second one consists in proving
and exploiting suitable Strichartz estimates (following
\cite{GSNT04,Miz08}) to get asymptotic stability.

The first step goes as follows: consider first the system with
$\epsilon=0$ and introduce action angle coordinates $(I,\alpha)$ for
the zero-th oscillator, thus one is reduced to a perturbation of a
Hamiltonian of the form
\begin{equation}
\label{h.lin}
\hs_0(I)+\sum_{k\not=0} \frac{p_k^2+q_k^2}{2}\ ,
\end{equation}
with $\hs_0$ a suitable function.  
If the perturbation does not contain terms linear in $(p,q)$ then the
manifold $p=q=0$ is invariant. So the idea is to iteratively eliminate
from the perturbation the terms linear in such variables. Furthermore
it also useful to eliminate the terms of order zero in $p,q$,which
depend on the angle $\alpha$ conjugated to $I$. This is expected to
be possible under the so called first Melnikov condition, namely
$$
\omega_0\not=1/n\ ,\quad \omega_0:=\frac{\partial\hs_0}{\partial
  I}\ ,\quad  n\in\Z\ .
$$ However we have not been able to find rigorous results on this
problem before the paper \cite{Gio12} in which Giorgilli proved the
convergence of the normal form in the case of Lyapunov periodic
orbits. The method by Giorgilli is based on his previous work
\cite{Gio01} (an improvement of Cherry's theorem \cite{Che}). Actually
it consists of a careful analysis (and estimate) of the formal
iterative procedure used to put the system in normal form, analysis
which allows to prove the convergence of such an iterative procedure.

Here we use a variant of Giorgilli's method. Theorem \ref{t.nor.for}
of the present paper differs from Giorgilli's one in the fact that we
are here in an infinite dimensional context and we also need here
to keep control of some weighted norms of the
perturbation. Furthermore, we have to study quite explicitly the
first two steps of the iterative procedure in order to have
a precise description of the linearization of the Hamiltonian at the
breather.

\vskip 10pt 

The dispersive step is more standard and consists of a variant of the
theory of \cite{KPS09}, which in turn is based on the previous results
\cite{SK05}, \cite{KKK}, \cite{PS08} (see also \cite{CucTar09}) and on
ideas by \cite{Miz08}. The only difference with respect to such works
rests in the fact that in our case the dispersion is of order
$\epsilon$ and we need to keep into account the dependence of all the
constants on $\epsilon$, thus we repeat, when needed some steps of the
proofs of such papers.

It is worth mentioning that the requirement of having a nonlinearity
starting with high degree is present also in all the quoted papers and
up to now there are no results on the case of analytic nonlinearities
with a potential vanishing at an order smaller then 8. It is probably
possible to weaken such a requirement by increasing the dimension of
the lattice. We also remark that the extension of the normal form
Theorem \ref{t.nor.for} to higher dimensions is straightforward, while
the adaptation of the dispersive part requires probably some
nontrivial work. Finally we point out that the theory of this paper
can also be adapted to deal with the model of \cite{Bam98} in which
the on site potential does not contain the quadratic term.

The paper is organized as follows: in Sect. \ref{state} we state the main
result; in Sect. \ref{pert} we state and prove the normal form result;
in Sect. \ref{de} we deal with the dispersive part of the proof and
conclude the proof of the main theorem; in Appendix \ref{technical}
we prove some technical lemmas needed in the part on normal form; in
Appendix \ref{techdis} we give some technical lemmas needed for the
dispersive part.

\vskip10pt
\noindent{\it Acknowledgments.} First I would like to thank
A. Giorgilli for some discussions on normal form theory and for
pointing to my attention his works. I also thank J. Villanueva and
H. Broer for some information on the normal form problem,
D. Pelinovski and A. Komech for some bibliographic indications on the
dispersive behavior of lattices.

\section{Statement of the main result}\label{state} 

We first introduce action angle variables
$$
(I,\alpha)\in\R_+\times\toro
$$ (here $\toro:=\R/\Z$ is the torus ) for the zero-th oscillator. We
recall that these variables are characterized by the following
properties: $I,\alpha$ are canonically conjugated, $\alpha$ is an
angle (i.e. $\alpha=\alpha+2\pi$), and the one particle Hamiltonian is
a function of $I$ only,
\begin{equation}
\label{2}
\frac{p_0^2+q_0^2}2+V(q_0)=\hs_0(I)\  ,
\end{equation}
with a suitable function $\hs_0$.  We recall that, if $V$ is analytic
and fulfills $V(q)=O(|q|^3)$ then also the variables $(I,\alpha)$ are
analytic in a domain of the form
$(0,C)\times\toro$, $C>0$, and then also
$\hs_0$ is analytic.

From now on we parametrize the phase space by using the coordinates
$(I,\alpha,p,q)$, $p=(p_k)_{k\not=0}$,
$q=(q_k)_{k\not=0}$. Furthermore we will use the collective notations
$\xi\equiv(p,q)$ and $\zeta\equiv(I,\alpha,\xi)$.

\vskip10pt

We denote by $\ell^{r}_s$ the space of the sequences $q\equiv (q_k)$
such that
\begin{equation*}
\norma{q}_{\ell^{r}_s}:=\left(\sum_{k}|q_k|^r\langle
k\rangle^{rs}\right)^{1/r} \ ,\quad s\in\R\ ,\quad 1\leq r<\infty
\end{equation*}
is finite. As usual $\langle k\rangle:=\sqrt{1+k^2}$, 
$\ell^{\infty}_s$ is defined by the sup norm.

We will also denote by $\bl^r_s:=\ell^r_s\oplus\ell^r_s$. If
$s=0$ we will write $\ell^{r}_0=:\ell^r$ and similarly for $\bl^r$.

We will use also spaces with exponential weights: we fix once for all a
positive $\beta$ and consider the spaces $\ell^{+}$, respectively
$\ell^-$ of the sequences such that the norm
\begin{align}
\label{norme}
\norma{ q }^2_+:=\sum_{k}\e^{\beta |k|}|q_k|^2\ ,\ {\rm
  respectively}\quad \norma{ q }^2_-:=\sum_{k}\e^{-\beta
  |k|}|q_k|^2\ ,
\end{align}
is finite. We will also denote $\Ph^{\pm}:=\ell^{\pm}\times\ell^{\pm}$.

\begin{remark}
\label{abuse}
We did not specify the range of the index $k$. Most of times it will
run over $\Z-\left\{0\right\}$, but sometimes over the whole
$\Z$. Every time this will be clear from the context.
{\it Furthermore, by abuse of notion we will say that a phase point
$\zeta\equiv(I,\alpha,\xi)\in\bl^r_s$ (or $\zeta\in\Ph^{\pm}$) if
$\xi\in\bl^r_s$ (or $\xi\in\Ph^{\pm}$).}
\end{remark}

Given two phase point $\zeta\equiv(I,\alpha,\xi)$ and
$\zeta'\equiv(I',\alpha',\xi')$ we define their distance according to
the different norms by
\begin{align*}
\d_{\bl^r}(\zeta;\zeta')&:=\max\left\{
\left|I-I'\right|;\left|\alpha-\alpha'\right|;
\norma{\xi-\xi'}_{\bl^{r}} \right\}
\\ \d_{\pm}(\zeta;\zeta')&:=\max\left\{
\left|I-I'\right|;\left|\alpha-\alpha'\right|; \norma{\xi-\xi'}_{\pm}
\right\}\ .
\end{align*}

Following \cite{KT} we say that a pair $(q,r)$ is admissible if $q\geq
6$, $r\geq 2$ and 
\begin{equation}
\label{a.p}
\frac{1}{q}+\frac{1}{3r}\leq \frac{1}{6}\ .
\end{equation}

All along the paper we will use the notation $a\sleq b$ to mean
``there exists a constant $C$, independent of all the relevant
quantities, such that $a\leq Cb$''. Sometimes, when needed or when
interesting, we will write explicitly the constant.

\vskip10pt

Denote by $b_0(I,t)$ the family of periodic solutions of the system
with $\epsilon=0$ defined by
\begin{equation}
\label{bzero}
b_0(I,t):=(I,\omega_0t+\alpha_0,0)=(I(t),\alpha(t),\xi(t))\ ,\quad
\omega_0:=\frac{\partial \hs_0}{\partial I}(I)\
\end{equation}
and by $\gamma_0:=\bigcup_{t}b_0(I,t)$ the corresponding trajectory,
then the main result of the paper is the following Theorem.

\begin{theorem}
\label{breather}
Assume that $V$ is analytic in a neighborhood of zero and that 
$V(q)=O(|q|^8)$ as $q\to0$, assume also that there exist positive
$C_{\omega_0}$, and $\Delta_1<\Delta_2$, such that the variables
$(I,\alpha)$ are real analytic in $[\Delta_1,\Delta_2]\times \toro$
and the following inequality holds 
\begin{equation}
\label{non.res}
\left|\omega_0(I)-\frac{1}{n}\right|\geq C_{\omega_0}\ ,\quad \forall
n\in\Z\ ,\quad\forall I\in[\Delta_1,\Delta_2]\ ,
\end{equation}
then there exists $\epsilon_*>0$, such that, for any
$0<\epsilon<\epsilon_*$ there exists a family of periodic solutions
$b_\epsilon(\I,t)$, $\I\in[\Delta_1,\Delta_2]$, of the system
\eqref{0.1}, with trajectories
$\gamma_\epsilon(\I):=\cup_{t}b_\epsilon(\I,t)$, having the following
properties:
\begin{itemize}
\item[i)] the distance between the unperturbed breather and the true
  breather is small:
  $\displaystyle{\d_+(\gamma_\epsilon(\I);\gamma_0(\I))}\sleq \sqrt
  \epsilon$,
\item[ii)]the family $\gamma_\epsilon(\I)$ is asymptotically
  stable. Precisely, fix $\delta>1/2$, then the following holds true:
  there exists $\epsilon_\delta>0$ such that, if
  $\epsilon<\epsilon_\delta$ and the initial datum $\zeta_0$ fulfills
\begin{equation}
\label{dato.i}
\inf _{\I\in[\Delta_1,\Delta_2]}\d_{\bl^2}(\zeta_0,\gamma_\epsilon(\I))
=:\mu<\epsilon^\delta\ ,
\end{equation}
then there exists an analytic function $\I(t)$ s.t.
\begin{itemize}
\item[ii.1)] for any admissible pair $(q,r)$ the function $t\mapsto
  \d_{\bl^r}(\zeta(t);\gamma_\epsilon(\I(t)))$ is of class $L^q_t$ and fulfills
\begin{equation}
\label{asybreather}
\norma{\d_{\bl^r}(\zeta(.);\gamma_\epsilon (\I(.)))}_{L^q_t}\sleq
\epsilon^{-1/q}\mu\ .  
\end{equation}
\item[ii.2)]$\displaystyle{|\I(t)-\I(0)|
  \sleq\frac{\mu^2}{\epsilon^{1/2}}}$ and
  $\displaystyle{\I_{\pm}:=\lim_{t\to\pm\infty}\I(t)}$ exists.
\end{itemize}
\end{itemize}
\end{theorem}

\section{Construction of the breather and normal form close to it}
\label{pert}

In this section $\Ph^{\pm}$ will always denote the space of the {\bf
  complex} sequences $\zeta=(I,\alpha,\xi)$
s.t. $\norma{\xi}_{\pm}<\infty$.

\subsection{Statement}\label{sta}

The main result of this section is the following theorem.

\begin{theorem}
\label{t.nor.for}
Assume that $V$ is analytic in a neighborhood of zero and
$V(q)=O(|q|^4)$ as $q\to0$, assume also that there exist positive
$C_{\omega_0}$, $\Delta_1<\Delta_2$, such that  such that the variables
$(I,\alpha)$ are real analytic in $[\Delta_1,\Delta_2]\times \toro$
and the following inequality holds 
\begin{equation}
\label{nonres.11}
\left|\omega_0(I)-\frac{1}{n}\right|\geq C_{\omega_0}\ ,\quad \forall
n\in\Z\ ,\quad\forall I\in[\Delta_1,\Delta_2]\ ,
\end{equation}
then there exists $\epsilon_0>0$, and $\forall |\epsilon|<\epsilon_0$
there exist complex neighborhoods $\U^{\pm}\subset\Ph^{\pm}$ of
$[\Delta_1,\Delta_2]\times \toro\times\{0\}$ and an analytic canonical
transformation $T:\U^{\pm}\to\Ph^{\pm}$ leaving invariant the space of
real sequences, with the following properties:
\begin{itemize}
\item[i)] there exists a positive $K_1$ s.t.
\begin{equation}
\label{dom.+}
\U^+\supset [\Delta_1-\frac{1}{K_1},\Delta_2+\frac{1}{K_1}]\times
\toro\times\left\{\xi\ :\ \norma\xi_{+}
\leq\frac{\sqrt\epsilon}{K_1} \right\}\ ,
\end{equation}
and
\begin{equation}
\label{dom.-}
\U^-\supset [\Delta_1-\frac{1}{K_1},\Delta_2+\frac{1}{K_1}]\times
\toro\times\left\{\xi\ :\ \norma\xi_{-}
\leq\frac{\sqrt\epsilon}{K_1} \right\}\ .
\end{equation}

\item[ii)] the transformed Hamiltonian $H\circ T$ has the
form
\begin{equation}
\label{H.tra}
H\circ T=\hs(I)+H_L+\V+\Zc\ ,
\end{equation}  
where 
\begin{itemize}
\item[ii.1)]
\begin{align}
\label{HL}
H_L:=\sum_{k\not=0}\frac{p_k^2+q_k^2}{2}+\epsilon\left[
  \sum_{k\not=-1,0}\frac{(q_{k+1}-q_k)^2}{2}+q_{-1}^2+q_1^2
  \right]\ ,
\\
\label{VC}
\Vc(q):=\sum_{k\not=0}V(q_k)\ ,
\end{align}
\item[ii.2)] $\hs(I)$ is an analytic function of $I$ fulfilling (with
  an $l$-dependent constant)
$$ \sup_{\U^-}\left|\frac{\partial^l (\hs-\hs_0)}{\partial
  I^l}\right|\sleq \sqrt{\epsilon}\ ,\quad\forall l\geq 0\ ;
$$
\item[ii.3)] $\Zc$ is such that its Hamiltonian vector field
$X\equiv(X_I,X_\alpha,X_{\xi})$ is analytic as a map
$X:\U^{-}\to\Ph^+$ and its components fulfill the following estimates
\begin{align}
\label{sti.i}
\sup_{(I,\alpha,\xi)\in\U^-}\left|X_I(I,\alpha,\xi)\right|\sleq
\epsilon^{1/2}\norma{\xi}_-^2 
\\
\label{sti.xi}
\sup_{(I,\alpha,\xi)\in\U^-}\norma{X_\xi(I,\alpha,\xi)}\sleq
\epsilon^{3/2}\norma{\xi}_- \ .
\end{align}

\end{itemize}
\item[iii)]  $T$ fulfills the estimates
\begin{align}
\sup_{(I,\alpha,\xi)\in\U^-}\left|T_I(I,\alpha,\xi)-I\right|\sleq
\epsilon ^{1/2}
\\
\label{sti.ti}
\sup_{(I,\alpha,\xi)\in\U^-}\left|T_\alpha(I,\alpha,\xi)-\alpha\right|\sleq
\epsilon^{1/2}
\\
\nonumber
\sup_{(I,\alpha,\xi)\in\U^-}\norma{T_\xi(I,\alpha,\xi)-\xi}_{+}\sleq
\epsilon 
\end{align}
where $T_I,T_\alpha,T_{\xi}$ are the different components of $T$.
\end{itemize}
\end{theorem}
\begin{remark}
\label{cor.bre}
The Hamiltonian $H\circ T$ admits the invariant manifold $\xi=0$ which
is foliated in periodic orbits. In the original coordinates such
periodic orbits are exponentially localized in space and in fact are
the breathers by MacKay and Aubry. Theorem \ref{t.nor.for} also
contains some information on the Hamiltonian close the breather,
information which is crucial for proving asymptotic stability.
\end{remark}

\begin{remark}
\label{lr.analit}
Since, for any $r\geq 1$, the embeddings 
$$
\Ph^+\hookrightarrow \bl^r\hookrightarrow\Ph^-
$$
are continuous, the transformation $T$ is analytic also as a map from
$\bl^r$ to itself.
\end{remark}

\subsection{Proof of theorem \ref{t.nor.for}}\label{proof.pert}

Before starting the construction it is useful to make the
following coordinate transformation:
\begin{align}
\label{trans.1}
z_k=\frac{p_k+\im {q_k}}{\sqrt{2}} 
\\
\nonumber
w_k=\frac{p_k-\im {q_k}}{\sqrt{2}} 
\end{align}
which transform the symplectic form to
$$
\di I\wedge \di \alpha+\im \sum_k \di z_k\wedge \di w_k\ .
$$ The transformation \eqref{trans.1} only multiplies the norms by a
constant, so it is enough to prove theorem \ref{t.nor.for} in the new
variables. In this section (and in appendix \ref{technical}) we will
denote by $\xi\equiv(z,w)$ the new complex variables; the collection
of the variables $(I,\alpha)$ will be denoted by $x\equiv(I,\alpha)$ .

In order to keep into account the different size of the different
variables we proceed as follows: fix some positive constants
$R_\alpha,R_I,$ and define $R_\xi:=\sqrt\epsilon$, then given a point
$\zeta\equiv (I,\alpha,\xi)$ we define its norms by
\begin{equation}
\label{norme2}
\snorma{\zeta}_{\pm}:=\max\left\{\frac{|I|}{R_I}, 
\frac{|\alpha|}{R_\alpha},\frac{\norma{\xi}_{\pm}}{R_\xi}
\right\}\ .
\end{equation}
Sometimes we will also denote 
\begin{equation}
\label{3.13}
\snorma{\xi}_{\pm}:=\frac{\norma\xi_{\pm}}{R_\xi}\ .
\end{equation}

The complex closed ball of radius $R$ and center $\zeta$ in such
topologies will be denoted by $B_{\pm}(R,\zeta)$.
\begin{remark}
\label{rem.1}
This is an $\epsilon$ dependent norm. The dependence of all the
constants on $\epsilon$ will be recorded, on the contrary the
quantities  $R_\alpha,R_I,$ will play no role and will be considered
as fixed. 
\end{remark}

We will develop perturbation theory in a complex neighborhood of the
domain 
\begin{equation}
\label{dom.g}
\G:=[\Delta_1,\Delta_2]\times \toro\times \{0\}\ni(I,\alpha,\xi)\ .
\end{equation}
We fix once for all a positive $R$. For
$\delta\in[0,1)$ we denote
\begin{equation}
\label{dom.g.1}
\G_{\delta}^{\pm}:=\bigcup_{\zeta\in\G}B_{\pm}(\zeta,R(1-\delta))\ .
\end{equation} 

We now define what we mean by normal form.

\begin{definition}
\label{def.nor.for}
For some $1>\delta\geq 0$, let $f=f(I,\alpha,\xi)$ be a Hamiltonian
function analytic on $\G^+_\delta$. The function $f$ will be said to
be in normal form if $f(I,\alpha,0)=0$ and $\di_\xi
f(I,\alpha,0)\equiv 0$, where $\di_\xi$ is the differential with
respect to $\xi$.
\end{definition}

\begin{remark}
\label{rem.2}
If a Hamiltonian function has the
form $$H=\hs(I)+\sum_{k\not=0}z_kw_k+ f$$ with $f$
in normal form then the manifold $\xi=0$ is invariant for
the dynamics and is foliated in periodic
orbits with frequency $\partial_I\hs(I)$. 
\end{remark}

To start with we introduce some notations. 

Given a Hamiltonian function $f=f(I,\alpha,\xi)$ we will denote 
\begin{align}
\label{0,1,2}
f^{(0)}(I,\alpha)&:=f(I,\alpha,0)\ ,\quad \left\langle
f^{(0)}\right\rangle(I):=
\frac{1}{2\pi}\int_0^{2\pi}f^{(0)}(I,\alpha)\di \alpha\ ,
\\
f^{(1)}(I,\alpha,\xi)&:=
\di_\xi f^{(1)}(I,\alpha,0,0)\xi
\\
&\equiv  \sum_{k\not=0}\left[\di_{z_k}
  f^{(1)}(I,\alpha,0,0)z_k+\di_{w_k} f^{(1)}(I,\alpha,0,0)w_k\right]
\ ,
\\
f^{(2)}&:=f-f^{(0)}-f^{(1)}\ ,
\end{align}
so that $f^{(2)}$ is in normal form.

Furthermore, for $f=f^{(1)}$ there exists a map
$f^1(I,\alpha)$ such that 
\begin{equation}
\label{p.1}
f^{(1)}(I,\alpha,\xi)=\left\langle f^1(I,\alpha);\xi \right\rangle:=
\sum_{k\not=0}\left(f^1_{z,k} z_k +f^1_{w,k} w_k   \right)\ ,
\end{equation}
where the scalar product is that of $\bl^2$.

Given a Hamiltonian function $\chi$ on $\G_\delta^{\pm}$ we will denote
by $X_\chi$ its Hamiltonian vector field, by
$\left[X_\chi\right]_\alpha\equiv \der \chi I$ its $\alpha$ component, and
similarly all the other components.

It is useful to introduce the operator $J$ (Poisson tensor) defined
by
$$
J\left(\begin{matrix}z\\w
\end{matrix}\right)
=
\left(\begin{matrix}
-\im w
\\
\im z
\end{matrix}\right)\ ,
$$
so that, with the above notations
$$
\left[X_{f^{(1)}}\right]_\xi\equiv Jf^1\ .
$$

To measure the size of the Hamiltonian vector fields of functions we
will use the following norms
\begin{align}
\label{normad}
\normad 
\delta
\chi&:=\frac{1}{R}max\left\{\sup_{\zeta\in\G^+_{\delta}}\snorma{X_\chi(\zeta)}_+,
\sup_{\zeta\in\G^-_{\delta}}\snorma{X_\chi(\zeta)}_-\right\} \ ,
\\
\label{normas}
\normas
{\delta} {\chi}&:=\frac{1}{R}\sup_{\zeta\in\G^-_{\delta}}\snorma{X_\chi(\zeta)}_+\ .
\end{align}

\begin{definition}
\label{functions}
A function whose Hamiltonian vector field is analytic as a map from
$\G_\delta^{\pm}$ to $\Ph^{\pm}$ will be said to be of class
$\A_\delta$. A function whose Hamiltonian vector field is analytic as
a map from $\G_\delta^{-}$ to $\Ph^{+}$ will be said to be of class
$\Sc_\delta$. 
\end{definition}

The key estimates which will be used in estimating the normal form are
given in the following lemma.
\begin{lemma}
\label{sti.varie}
Let $f\in\Sc_d$ and $g\in\A_{d}$ be analytic functions, then one has 
\begin{align}
\label{4.3.1}
\normas{d+d_1}{\poisson f{g}}\leq\frac2{d_1}\normas
d f \normad d{g}\ .
\\
\label{sti.parti.21}
\normas d{f^{(0)}}\leq \normas d{f}\ ,\quad 
\normas d{f^{(1)}}\leq \normas d{f}\ ,
\quad
\\
\label{stime.parti.111}
\normas {d+d_1}{f^{(2)}}\leq \frac{1}{d_1^2}\normas d{f}\ ,
\\
\label{poi.part.1.1}
\normas d{\poisson{f^{(1)}}{g^{(2)}}^{(1)}}\leq
\frac{3}{1-d}\normas d{f^{(1)}} \normad d{g^{(2)}} \ .
\end{align}
\end{lemma}
The proof will be given in Appendix \ref{technical}.

In particular the estimate \eqref{poi.part.1.1} in which there is no
$d$ at the denominator (but $1-d$, which is bounded away from zero) is
the key for the convergence of the normal form procedure. 

The canonical transformation increasing by one the order of the non
normalized part of the Hamiltonian will be constructed as the Lie
transform generated by auxiliary Hamiltonian functions of the form
$\chi^{(1)}+\chi^{(0)}\in\Sc_d$ with some $d$. 

\begin{remark}
\label{flow1}
Denote by $\Phi_\chi^t$ the flow of the Hamiltonian vector field of a
Hamiltonian function $\chi\in\Sc_d$, then by standard existence and uniqueness
theory one has that, if $\normas d\chi\leq d_1 $ with $0\leq
d_1<1-d$ then $\Phi^t_\chi$ exists at least up to time 1,
furthermore one has
\begin{equation}
\label{esttchi}
\sup_{\zeta\in\G^-_{d+d_1}}\snorma{\Phi^t_{\chi}(\zeta)-\zeta}_+\leq
R\normas d\chi \ .
\end{equation}
\end{remark}
\begin{remark}
\label{taylor}
By standard Hamiltonian theory, for any smooth $f$ one has
$$
\frac{\di}{\di t}f\circ \Phi^t_{\chi}=\poisson \chi f\circ
\Phi^t_\chi\ ,
$$
thus, defining the sequence $f_{(l)}$ by
$$
f_{(0)}:=f\ ,\quad f_{(l)}:=\poisson \chi{f_{(l-1)}}\ ,\quad l\geq 1\ ,
$$
one has, for any $N\geq 0$
\begin{equation}
\label{restotaylor}
f\circ
\phi^t_\chi=\sum_{l=0}^{N}\frac{f_{(l)}}{l!}+\int_0^1\frac{(1-s)^{N}}{N!}
f_{(N+1)}\circ\Phi^s_{\chi}\, \di s\ .
\end{equation}
\end{remark}

\begin{lemma}
\label{l.2}Let
 $\chi\in\Sc_d$, with $0\leq d_1<1-d$
 and let $f\in\A_ d$; fix $0<
d_1<(1- d)$, assume $\normas d {\chi}\leq d_1/3$, then one has 
\begin{align}
\label{rest.tayl}
\normas
    {d+d_1}{f\circ\Phi^1_\chi-\sum_{j=0}^{N}\frac{1}{j!}f_{(j)}}\sleq
    \frac{1}{d_1^{N+1}}\normad d f(\normas d\chi)^{N+1}\ .
\end{align}
\end{lemma}
The proof will be given in Appendix \ref{technical}.

As usual, in order to find the generating function for the normalizing
transformation one has to solve a cohomological equation, which in our
case will have the form 
\begin{equation}
\label{homo.1}
\poisson{H_{lin}}{\chi}= \Psi^{(0)}+\Psi^{(1)}\ , 
\end{equation}
where $\chi$ is the unknown, $\Psi^{(0)}$, $\Psi^{(1)}$ are given functions,
\begin{equation}
\label{hlin}
H_{lin}(I,\xi):=\hs(I)+\sum_{k\not=0} z_kw_k
\end{equation}
and $\hs$ is a function of the action $I$ only. The last estimate we
need before starting the recursive construction of the normal form is
contained in the following lemma, which will proved in section
\ref{technical}.

\begin{lemma}
\label{sol.homo}
Let $1>\delta\geq0$ be given and consider the equation
\eqref{homo.1}. Assume that $\langle\Psi^{(0)}\rangle=0$, that
$\Psi^{(0)}\in\Sc_\delta$, $\Psi^{(1)}\in\Sc_\delta$ and
$\hs\in\Sc_\delta$. Denote $\omega(I):=\der {\hs}I$, and assume
that on $\G^-_{\delta}$ one has
\begin{equation}
\label{nonres.1}
\left|\omega(I)-\frac{1}{n}\right|\geq C_\omega>0\ ,\quad \forall
n\in\Z\ ,
\end{equation}
then the cohomological equation \eqref{homo.1} has a solution
$\chi=\chi^{(0)}+\chi^{(1)}$ which fulfills
\begin{align}
\label{sti.homo}
\normas \delta{\chi^{(0)}}&\leq \pi \left[\sup_{\G^{-}_\delta}
  \frac{1}{\omega(I)} \right] \normas \delta{\Psi^{(0)}}\sleq
\normas\delta{\Psi^{(0)}} \ , \\
\label{sti.homo.1}
 \normas \delta{\chi^{(1)}}&\leq 2\pi \left[\sup_{\G^{-}_\delta}
   \frac{1}{\omega(I)\left|\sin\left(\frac{\pi}{\omega(I)}
     \right)\right|} \right] \normas \delta{\Psi^{(1)}}\sleq
 \normas\delta{\Psi^{(1)}} \ .
\end{align}
\end{lemma} 

We now check the analyticity properties of the vector field of the
Hamiltonian \eqref{0.1}, which we rewrite as 
\begin{align}
\label{0.1.1}
\Hc_0&:=H=H_0+\Zc_2+\resto^{(1)}_0+\resto^{(0)}_0
\\
\label{0.1.h0}
H_0&:=\hs_0(I)+\sum_{k\not=0}z_kw_k\ ,
\\
\label{0.1.Z}
\Zc_2&:=\epsilon\left[
  \sum_{k\not=-1,0}\frac{(q_{k+1}-q_k)^2}{2}+q_{-1}^2+q_1^2 \right]
+\sum_{k\not=0}V(q_k) 
\\
\label{0.1.r0}
\resto_0^{(1)}&:=-\epsilon q_0(I,\alpha)\left[q_{-1}+q_1\right]\ ,\quad
\resto_0^{(0)}:=\epsilon [q_0(I,\alpha)]^2\ ,
\end{align}
where $q_k=q_k(z_k,w_k)$ for $k\not=0$. 

\begin{lemma}
\label{primo.0}
There exists 
$\epsilon_{*1}>0$ such that, if $|\epsilon|<\epsilon_{*1}$, then
$\Zc_2\in\A_{0}$, while $\resto_0^{(1)},\resto_0^{(0)}, \hs_0(I)\in\Sc_0$ and
the following estimates hold
\begin{align}
\label{sti.0.00}
\normad0{\Zc_2}\sleq \epsilon\ ,\quad \normas0{\resto^{(0)}_0}\sleq
\epsilon\ ,
\quad
\normas0{\resto^{(1)}_0}\sleq \sqrt\epsilon\ .
\end{align}
\end{lemma}
The very simple proof is left to the reader.

We proceed in constructing the canonical transformation putting the
system in normal form.  To this end we have to fix a sequence of
domains in which the transformed Hamiltonians will be defined. Thus
fix
\begin{equation}
\label{delta}
\delta_j:=\delta \e^{-j}\ ,\quad \tilde \delta_r:=\sum_{j=1}^{r}\delta_j
\ ,\quad \delta:=(\e-1)/2\ ,
\end{equation} 
so that $\sum_{j\geq1}\delta_j=1/2$.

The first two steps of the normalizing procedure have to be performed
in detail in order to keep the needed information on the linearization
of the equations at the breather.

\begin{lemma}
\label{primi2}
Assume that \eqref{nonres.1} holds, then there exist positive
$\epsilon_{*2}$ such that, for any $|\epsilon|<\epsilon_{*2}$, there
exists an analytic canonical transformation $T_2:\G^-_{\tilde \delta_2} \to
\G^-_{0}$ which restricts to an analytic transformation
$T_2:\G^+_{\tilde \delta_2} \to \G^+_{0}$ such that
\begin{align}
\label{hc2}
\Hc_2:=\Hc_0\circ T_2=H_2+\Zc_2+\resto_2
\\
\nonumber
H_2=\hs_2(I)+\sum_{k\not=0}z_kw_k \ ,\quad \hs_2=\hs_0+h_2 
\end{align}
and the following estimates hold
\begin{equation}
\label{esti.prime2}
\normas{\tilde \delta_2}{\resto_2}\sleq\epsilon^{3/2}\ ,\quad
\normas{\tilde\delta_2}{h_2}\sleq \epsilon  
\end{equation}
Furthermore one has
$T_2=(\uno+\Tr_1)(\uno+\Tr_2)$ with $\Tr_j:\G^-_{\tilde
  \delta_j}\to\G^+_{\tilde \delta_{j-1}} $
($j=1,2$) analytic and fulfilling
\begin{equation}
\label{t.1.1}
\sup_{\zeta\in \G^-_{\tilde \delta_j}}\snorma{\Tr_j(\zeta)}_+\sleq
(\sqrt\epsilon)^j
\ .
\end{equation} 
\end{lemma}
\begin{remark}
\label{r.primi2}
The important fact is that, up to order $\epsilon^{3/2}$ there are no
contributions correcting $\Zc_2$, whose form is explicitly known.
\end{remark}

\proof We proceed in two steps, each one increasing by $\epsilon^{1/2}$
the order of the non normalized part of the Hamiltonian.

Let $\chi_{1}^{(1)}$ be the solution of the cohomological equation 
$$
\left\{H_0;\chi^{(1)}_{1}\right\}= \resto_0^{(1)}\ ,
$$
so that $\normas0{\chi^{(1)}_{1}}\sleq \sqrt\epsilon$. Use
$\Phi^1_{\chi^{(1)}_{1}}$ to transform the Hamiltonian, then one has
\begin{align*}
\Hc_{01}:=\Hc_0\circ \Phi^1_{\chi^{(1)}_{1}}&=
H_0+\left\{\chi^{(1)}_{1};H_0 \right\} +\frac{1}{2}
\left\{\chi^{(1)}_{1};\left\{\chi^{(1)}_{1};H_0 \right\}\right\} \\ &+
\frac{1}{2}
\int_0^1(1-s)^2\left\{\chi^{(1)}_{1};\left\{\chi^{(1)}_{1};\left\{\chi^{(1)}_{1};H_0
\right\}\right\} \right\} \circ \Phi^s_{\chi^{(1)}_{1}}\di s \\ &+
\Zc_2+\int_0^1\left\{\chi^{(1)}_{1};\Zc_2\right\}\circ\Phi^s_{\chi^{(1)}_{1}}
\di s \\ &+ \resto_0^{(0)}+
\int_0^1\left\{\chi^{(1)}_{1};\resto_0^{(0)}\right\}\circ
\Phi^s_{\chi^{(1)}_{1}} \di s \\ &+
\resto^{(1)}_0+\left\{\chi^{(1)}_{1};\resto^{(1)}_0
\right\}+\int_0^1\left\{\chi^{(1)}_{1};\left\{\chi^{(1)}_{1};\resto^{(1)}_0
\right\}\right\}\circ \Phi^s_{\chi^{(1)}_{1}}\di s \\ &=H_0+\frac{1}{2}
\left\{\chi^{(1)}_{1};\resto^{(1)}_0\right\} \\ &+
\int_0^1\left[1-\frac{(1-s)^2}2\right]\left\{\chi^{(1)}_{1};\left\{\chi^{(1)}_{1};\resto^{(1)}_0\right\}
\right\} \circ \Phi^s_{\chi^{(1)}_{1}}\di s \\ &+
\Zc_2+\int_0^1\left\{\chi^{(1)}_{1};\Zc_2\right\}\circ\Phi^s_{\chi^{(1)}_{1}}
\di s \\ &+ \resto_0^{(0)}+
\int_0^1\left\{\chi^{(1)}_{1};\resto_0^{(0)}\right\}\circ
\Phi^s_{\chi^{(1)}_{1}} \di s
\end{align*}
By using \eqref{4.3.1} and lemma \ref{tra} the seventh line, the
integral at the eighth line and the integral at the ninth line have a
norm which is estimated by a constant times $\epsilon^{3/2}$. It
remains to estimate the Poisson bracket
$$
\left\{\chi^{(1)}_{1};\resto^{(1)}_0\right\}\equiv
\left\{\chi^{(1)}_{1};\resto^{(1)}_0
\right\}^{(0)}+\left\{\chi^{(1)}_{1};\resto^{(1)}_0\right\}^{(2)} \ .
$$ By equations \eqref{4.3.1} and \eqref{sti.parti.21} the first term at
r.h.s. has norm $\normas..$ of order $\epsilon$ and is independent of
$\xi$. We are now going to prove that

\begin{equation}
\label{sti.second}
\normas{\delta_1}{\left\{\chi^{(1)}_{1};\resto^{(1)}_0\right\}^{(2)}}\sleq
\epsilon^2\ .
\end{equation}
Denote $f:=\left\{\chi^{(1)}_1;\resto^{(1)}_0\right\}^{(2)}\ ;
$ using the further notation
$$
\chi^{(1)}_{1}=\left\langle \chi^1;\xi\right\rangle\ ,\quad
\resto_0^{(1)} =\langle\resto_0^{1};\xi\rangle\ ,
$$
one has
\begin{equation}
\label{f.r.4}
f=\left\langle \der{\chi^1}I;\xi \right\rangle\left\langle
\der{\resto^1_0}\alpha;\xi \right\rangle- 
\left\langle \der{\chi^1}\alpha;\xi \right\rangle\left\langle
\der{\resto^1_0}I;\xi \right\rangle\ ,
\end{equation}
so that 
\begin{equation}
\label{f.r.5}
-\left[X_f\right]_I= \left\langle \frac{\partial^2 \chi^1 }{\partial
  \alpha\partial I};\xi \right\rangle\left\langle
\der{\resto^1_0}\alpha;\xi \right\rangle+\text{similar\ terms}\ .
\end{equation}
By the definition of $\normas..$, one has 
$$
\frac{1}{RR_\xi}\sup_{\G^-_0}\norma{\chi^1(I,\alpha)}_+\leq
\normas0{\chi^{(1)}_{1}} \ ,
$$ 
and therefore on $\G^-_{\delta_1}$, by Cauchy estimate, one has 
$$
\frac{1}{R R_\xi}\norma{\frac{\partial^2 \chi^1 }{\partial
  \alpha\partial I}   }_+\leq
\frac{2}{R^2R_IR_\alpha\delta_1^2}\normas0{\chi^{(1)}_{1} }\ , 
$$
which gives 
$$
\frac{1}{RR_\alpha }\norma{\frac{\partial^2 \chi^1 }{\partial
  \alpha\partial I}   }_+\leq
\frac{2}{R^2R_IR_\alpha^2\delta_1^2}\normas0{\chi^{(1)}_{1} }R_\xi\ .
$$
Inserting in \eqref{f.r.5} and taking into account that
$R_\xi=\sqrt\epsilon$ one has
\begin{equation}
\label{sti.1.2}
\sup_{\G^-_{\delta_1}}\frac{1}{RR_\alpha}\left|\left\langle
\frac{\partial^2 \chi^1 }{\partial \alpha\partial I};\xi \right\rangle
\right| \leq \frac{2}{R^2R_IR_\alpha^2\delta_1^2}\normas0{\chi^{(1)}_{1}
}R^2_\xi\sleq \epsilon\sqrt\epsilon \ .
\end{equation}
In a similar way one gets
$$
\left|\left\langle
\der{\resto^{1}_0}\alpha;\xi \right\rangle\right|\sleq \epsilon^{3/2}\ .
$$ from which one has that the first term of \eqref{f.r.5} is of order
$\epsilon^3$.  All the other terms can be estimated similarly getting
the wanted estimate for the $\alpha$ and the $I$ components of the
vector field.

Concerning the $\xi$ component of the
vector field one has
\begin{equation}
\label{sti.ancora}
\left[X_f\right]_\xi=J\der {\chi^1}I \left\langle
\der{\resto^{1}_0}I;\xi \right\rangle+\text{similar\ terms}\ . 
\end{equation}
In particular, acting as above one immediately proves that on
$\G^-_{\delta_1}$, 
$$
\left|\left\langle
\der{\resto^{1}_0}I;\xi \right\rangle\right|\sleq \epsilon^{3/2}\ .
$$
We now add the estimate of the derivative of $\chi^1$:
$$
\sup_{\G^-_{\delta_1}}\norma{J\frac{\partial \chi^1}{\partial
    I}}_+\leq
\frac{1}{RR_I\delta_1}\sup_{\G^-_{0}}\norma{J\chi^1}_+
=\frac{1}{RR_I\delta_1}\sup_{\G^-_{0}}\norma{[X_{\chi^{(1)}} ]_\xi}_+
\sleq \epsilon^{1/2}R_\xi\ ,
$$
dividing by $RR_\xi$, one gets that the norm $\normas{\delta_1}.$ of
the first term of \eqref{sti.ancora} is of order
$\epsilon^2$. Considering all the other terms one gets \eqref{sti.second}.

We have thus shown that after this transformation the Hamiltonian has
the form
\begin{equation}
\label{hc.1}
\Hc_{1}=H_0+\Zc_2+\resto_{1}+\tilde\resto_2
\end{equation}
where 
\begin{equation}
\label{resto1}
\resto_{1}\equiv\resto_{1}^{(0)}:=\resto^{(0)}_0
+\left\{\chi^{(1)}_{1};\resto^{(1)}_0 \right\}^{(0)}\ ,\quad
\normas{\delta_1}{\tilde \resto_2}\sleq \epsilon^{3/2}
\end{equation}
We now perform the second step removing the part of $\resto_1$
dependent on $\alpha$. To this end define 
$$
\Psi_2\equiv\Psi_2^{(0)}:=\resto_1-\langle\resto_1\rangle\ ,
$$ and define $\chi_2\equiv \chi^{(0)}_2$ as the solution of the
cohomological equation
$\left\{H_0;\chi_2\right\}=\Psi_2$. Transforming $\Hc_1$ one gets
\begin{align*}
\Hc_2&:=\Hc_1\circ\Phi^1_{\chi_2}
=H_0+\int_0^1s\left\{\chi_2;\Psi_2\right\}\circ \Phi^s_{\chi_2}\di s
\\ &+\langle\resto_1\rangle+\Zc_2+
\int_0^1\left\{\chi_2;\Zc_2\right\}\circ \Phi^s_{\chi_2}\di s
\\ &+\int_0^1\left\{\chi_2;\langle\resto_1\rangle\right\}\circ \Phi^s_{\chi_2}\di s+
\tilde\resto_2\circ\Phi^1_{\chi_2} \ .
\end{align*}
Defining $h_2:= \langle\resto_1\rangle$, and $\resto_2$ to be the sum
of the various integrals and of $\tilde\resto_2\circ\Phi^1_{\chi_2}$
and estimating the different terms, one immediately gets the thesis.
\qed

We are now ready to state the iterative lemma which is the heart of
the proof.

\begin{lemma}
\label{iter}
(Iterative Lemma). Assume that on $[\Delta_1,\Delta_2]$ the non
resonance condition \eqref{nonres.11} holds, then there exist positive
constants $\epsilon_*,C_1,C_2,C_3,C_4,K$, such that the following
holds true: for any $r\geq 2$ and any $\epsilon$ with
$|\epsilon|<\epsilon_*$ there exists a canonical transformation
$T_r:\G^{-}_{\tilde \delta_{r}}\to \G^{-}_{\tilde \delta_{r-1}}$,
which restricts to an analytic transformation $T_r:\G^{+}_{\tilde
  \delta_{r}}\to \G^{+}_{\tilde \delta_{r-1}}$ s.t.
\begin{equation}
\label{Hcr}
\Hc_r:=\Hc_0\circ T_r=H_r+\Zc_r+\resto_r\ ,
\end{equation}
where
\begin{align}
\label{Hr}
H_r:=\hs_r(I)+\sum_{k\not=0}z_kw_k\ ,\quad
\hs_r:=\hs_0+h_1+...+h_r\ ,
\\
\label{zcr}
\Zc_r=\Zc_2+Z_3+...+Z_r\ ,
\\
\label{tr}
T_r=(\uno+\Tr_1)\circ...\circ(\uno+\Tr_r)
\end{align}
and $Z_j$ is in normal form for all $j$'s. The following
estimates hold
\begin{align}
\label{restor}
\normas{\tilde\delta_r}{\resto_r}\leq C_1(K\sqrt\epsilon)^{r+1}
\\
\label{hsj}
\normas{\tilde\delta_j}{h_j}\leq C_1(K\sqrt\epsilon)^{j}
\\
\label{zj}
\normas{\tilde\delta_j}{Z_j}\leq
\frac{C_2}{\delta_j^3}(K\sqrt\epsilon)^{j}
\\
\label{tj}
\sup_{\zeta\in\G^-_{\tilde\delta_j}}\snorma{\Tr_j(\zeta)}_+\leq C_3
(K\sqrt\epsilon)^{j} \ .
\end{align}
Furthermore one has
\begin{align}
\label{hsr}
\normas{\tilde\delta_r}{\hs_r-\hs_0}\leq C_4\epsilon 
\\
\label{zcr1}
\normad{\tilde\delta_r}{\Zc_r}\leq C_4\epsilon<1
\end{align}
\end{lemma}
\proof For $r=2$ the lemma coincides with lemma \ref{primi2}. We
assume it is true for some $r$ and we prove it for $r+1$. 

Define 
\begin{align}
\label{tiz}
h_{r+1}:=\langle\resto_{r}^{(0)}\rangle\ ,
\\
\label{psir}
\Psi_{r+1}:=\resto^{(0)}_r-\langle\resto_{r}^{(0)}\rangle+\resto_{r}^{(1)}
\ .
\end{align}
So in particular $h_{r+1}$ satisfies \eqref{hsj}, and one has
\begin{align}
\label{rec.1}
\normas{\tilde\delta_r}{\Psi^{(1)}_{r+1}}\leq
\nr\ ,\quad
\normas{\tilde\delta_r}{\Psi^{(0)}_{r+1}}\leq 2\nr\ ,
\end{align} 

Then we define $\chi_{r+1}$ to be the solution of 
\begin{equation}
\label{homo.r}
\left\{ H_r,\chi_{r+1}\right\}=\Psi_{r+1}\ ;
\end{equation} 
remark that, by \eqref{hsr}, provided $\epsilon$ is small enough
(uniformly in $r$), $\omega_r:=\der{\hs_r}I$ satisfies
\eqref{nonres.1} with a smaller constant $C_{\omega}$ independent of
$r$.  Therefore $\chi_r$ exists and fulfills
\begin{equation}
\label{chir1}
\normas{\tilde\delta_r}{\chi_{r+1}^{(0)}}\sleq
\nr \ ,\quad
\normas{\tilde\delta_r}{\chi_{r+1}^{(1)}}\sleq \nr\ ,
\end{equation}
with constants which are independent of $r$ (as all the constants that
will be suppressed using the symbol $\sleq$). 

We define now $\Tr_{r+1}:=\Phi^{1}_{\chi_{r+1}}-\uno$ and 
\begin{align}
\label{zr+1}
&Z_{r+1}:=\resto^{(2)}_r+\left\{\chi_{r+1}^{(0)};\Zc_r\right\}
+\left\{\chi_{r+1}^{(1)} ;\Zc_r
\right\} ^{(2)}
\\
\label{rer+1}
&\resto_{r+1}:= \left\{\chi_{r+1}^{(1)} ;\Zc_r
\right\} ^{(1)}
\\
\nonumber
&+\int_0^1\left[\left\{\chi_{r+1};\resto_r\right\}+
 (1-s)\left[
    \left\{\chi_{r+1};\left\{\chi_{r+1};\Zc_r\right\}\right\}-
\left\{\chi_{r+1};\Psi_{r+1} 
  \right\}
 \right] \right]\circ\Phi^s_{\chi_{r+1}}\di s 
\end{align}
so that $\Hc_{r+1}$ has the wanted form. We are now going to estimate
the different terms in order to prove that the estimates
\eqref{restor}-\eqref{zcr1} hold at level $r+1$.

Concerning $Z_{r+1}$ we estimate the last term, which is the worst one:
\begin{align}
\label{es.ric.l}
\normas{\tilde\delta_{r+1}}{\left\{\chi_{r+1}^{(1)} ;\Zc_r
\right\} ^{(2)}} \leq \frac{4}{\delta_{r+1}^2}\normad{\tilde\delta_{r}+\delta_{r+1}/2}{\left\{\chi_{r+1}^{(1)} ;\Zc_r
\right\}}
\\
\nonumber
\leq \frac{4}{\delta_{r+1}^2}
\frac{2}{\delta_{r+1}}
\normas{\tilde\delta_r}{\chi_{r+1}} \normad{\tilde\delta_r}{\Zc_r}
\sleq  \frac{1}{\delta_{r+1}^3} \normad{\tilde\delta_r}{\Zc_r}\nr
\end{align}

Adding the other estimates one gets 
\begin{equation}
\label{zr+11}
\normas{\tilde\delta_{r+1}}{Z_{r+1}}\sleq
\left(\frac{1}{\delta_{r+1}^2}+  \frac{\nz}{\delta_{r+1}}+
\frac{\nz}{\delta^3_{r+1}} 
\right) \nr\sleq \frac{1}{\delta_{r+1}^3}{\nr}\ ,
\end{equation}
which, provided one chooses $C_2$ to be $C_1$ times the constant not
written in the last of \eqref{zr+11} gives \eqref{zj} at level $r+1$.

We come to $\resto_{r+1}$. All the terms can be estimated in a
straightforward way using lemmas \ref{sti.varie}, \ref{l.2} and
\ref{tra} giving, 
\begin{align}
\label{resto.dif1}
\normas{\tilde \delta_{r+1}}{\resto_{r+1}}\sleq
\nr\left[\nz+\frac{\nr}{\delta_{r+1}}+\frac{\nr}{\delta_{r+1}^2}\nz
  \right]\ ,
\end{align}
Calling $\cg$ the constant making true \eqref{resto.dif1} one has  
\begin{align}
\label{resto.dif}
\normas{\tilde \delta_{r+1}}{\resto_{r+1}}\leq \cg
C_1(K\sqrt\epsilon)^{r+1} \left[
  C_4\epsilon+\frac{C_1(K\sqrt\epsilon)^{r+1}}{\delta
    \e^{-(r+1)}}+\frac{C_1(K\sqrt\epsilon)^{r+1}}{\delta^2
    \e^{-2(r+1)}}C_4\epsilon  \right] \ .
\end{align}
Taking $\epsilon_*$ small enough one can make the square bracket
smaller then $(C_4+C_1)\epsilon$, which shows that \eqref{restor}
is fulfilled at order $r+1$ if one defines $K:=\cg (C_4+C_1)$. Remark
that actually one increases the order of the perturbation by
$\epsilon$ at every step, however we made the choice of estimating
$\epsilon$ by $\sqrt\epsilon$ in order to be able to give a
formulation of the theorem which is also valid in the case $r=1,2$. 

All the other estimates are simpler and are omitted.

The key point in getting the estimate \eqref{restor}, which in turn is
the key to get the convergence, rests in the fact that we separated
from $\resto_{r+1}$ the second two terms of \eqref{zr+1} and
furthermore in the fact that the first term of \eqref{rer+1} fulfills
the improved estimate \eqref{poi.part.1.1}.
\qed 

\vskip20pt
\noindent{\it Proof of Theorem \ref{t.nor.for}.} First remark that,
due to the uniformity of the estimates, in Lemma \ref{iter} one can
pass to the limit $r\to\infty$, getting a transformation
$T:=T_{\infty}$ which is defined on
$\G^{\pm}_{3/8}\subset\G^{\pm}_{1/2}$, which puts the system in normal
form. 

To get the estimate \eqref{sti.ti} remark first that, from \eqref{tj},
one has
\begin{equation}
\label{sti.ti.infy}
\snorma{\uno-T}_{+}\sleq \sum_{j\geq 1} (K\sqrt\epsilon)^{j}\sleq
\epsilon^{1/2} \ ,
\end{equation} 
then \eqref{sti.ti} is just a component wise formulation of
\eqref{sti.ti.infy}.

To estimate $X$, first remark that $\Zc:=\lim_{r\to\infty}\Zc_r$ is
defined and analytic in $\G^-_{3/8}$ and fulfills
$$
\normas{\frac{3}{8}}{\Zc}\sleq \sum_{j\geq
  3}\normas{\frac{3}{8}}{\Zc_j}\sleq \epsilon^{3/2}\ ,
$$
which, written component wise, gives
\begin{align}
\label{pr.n.1}
\sup_{\zeta\in\G_{\frac{3}{8}}^-}\left|X_I(\zeta)\right|\sleq
\epsilon^{3/2}\ ,
\\
\label{pr.n.2}
\sup_{\zeta\in\G_{\frac{3}{8}}^-}\left|X_\alpha(\zeta)\right|\sleq
\epsilon^{3/2}\ ,
\\
\label{pr.n.3}
\sup_{\zeta\in\G_{\frac{3}{8}}^-}\norma{X_\xi(\zeta)}_+\sleq
\epsilon^{2}\ .
\end{align}

Since $\Zc$ is in normal form one has $X_I(I,\alpha,0)=\d
_{\xi}X_I(I,\alpha,0)=0$, and thus, using the standard formula for the
remainder of the Taylor expansion, one has
$$
X_I(I,\alpha,\xi)=\int_0^1(1-\tau)\d^2_\xi X_I(I,\alpha,\tau\xi)(\xi,\xi)\d
\tau \ ;
$$
Using the analyticity of $X_I$ as a function of $\xi$ in the domain
$\norma{\xi}_{-}\leq 3\epsilon^{1/2}/8$ one gets
$$
\sup_{\zeta\in\G^-_{\frac{3}{8}}\cap \{\norma{\xi}_-\leq
  \frac{\sqrt\epsilon}{4}\}} \norma{\d^2_\xi X_I(\zeta)} \leq
\frac{2}{(\sqrt{\epsilon}/4)^2}
\sup_{\zeta\in\G_{\frac{3}{8}}}\left|X_I(\zeta)\right| \sleq
\sqrt\epsilon\ ,
$$
where the norm at the first term is for $\d^2_\xi X_I$ considered as
quadratic form on the space of the $\xi$'s endowed by the norm
$\norma._-$. This proves \eqref{sti.i}. 

Similarly, using 
$$
X_\xi(I,\alpha,\xi)=\int_0^1\d_\xi X_\xi(I,\alpha,\tau\xi)\xi\d \tau
$$
equation \eqref{pr.n.3} and Cauchy estimate for the differential, one
gets \eqref{sti.xi}.\qed

\section{Dispersive estimates}\label{de}

We first establish decay and Strichartz estimates
for the group generated by the linear operator representing the first
order normal form, namely for the flow of the linear system with
Hamiltonian 
\begin{equation}
\label{d.h.l}
H_L:=\sum_{k\not=0}\frac{p_k^2+q_k^2}{2}+\frac{\epsilon}{2}\sum_{k\not=-1,0}(q_{k+1}-q_k)^2+\epsilon(q_1^2+q_{-1}^2)\ .
\end{equation}

\subsection{Linear local decay estimates}\label{le}

In order to prove the decay estimates it is useful to remark that the
system \eqref{d.h.l} consists of two decoupled systems, the first
consisting of the left hand part of the chain and having Hamiltonian
\begin{equation}
\label{d.h.l1}
H_{L1}:=\sum_{k\leq
  -1}\frac{p_k^2+q_k^2}{2}+\frac{\epsilon}{2}\sum_{k=-\infty}^{-1}
(q_{k+1}-q_k)^2\ ,
\end{equation}
where the phase space variables are $(p_k,q_k)_{k\leq -1}$, while $q_0\equiv
0$. Analogously the second system consists of the right hand part of the
chain. Furthermore the system \eqref{d.h.l1} can be viewed as the
restriction of the system with Hamiltonian
\begin{equation}
\label{d.h.l2}
H_S:=\sum_{k\in\Z}\frac{p_k^2+q_k^2}{2}+\frac{\epsilon}{2}\sum_{k\in\Z}(q_{k+1}-q_k)^2
\end{equation}
to skew symmetric sequences, namely the space of the sequences
$(p_k,q_k)_{k\in\Z}$ such that $p_k=-p_{-k}$, $q_k=-q_{-k}$. The same
is true for the system describing the right hand part of the chain. 

Thus we start by establishing the needed decay estimates for the
restriction of \eqref{d.h.l2} to skew-symmetric sequences (actually
when needed we will explicitly assume skew-symmetry of the
sequences). 

The system \eqref{d.h.l2} is a Klein Gordon chain with small
dispersion, so we actually follow the procedure of \cite{SK05} and
\cite{KKK} just keeping into account that we need estimates uniform
in $\epsilon$ and that we are just interested in skew-symmetric
sequences. 

 Consider the Hamilton equations of \eqref{d.h.l2}, namely
\begin{align}
\label{h.l.eq}
\dot p_k&=-q_{k}+\epsilon(q_{k+1}+q_{k-1}-2q_k)
\\
\nonumber
\dot q_k&=p_k
\end{align}  
and denote by $S^0_\epsilon(t)$ its evolution operator, namely the operator that
to $\xi\equiv (p,q)$ associates the value at time $t$ of the solution
with initial datum $\xi$.

First remark that, by conservation of energy one has that, for
$\epsilon$ small enough the system \eqref{h.l.eq} is globally well
posed in $\bl^2$ and the inequality
\begin{equation}
\label{c.l2}
\norma{S^0_\epsilon(t)\xi}_{\bl^2}\sleq \norma{\xi}_{\bl^2}\ 
\end{equation}
holds.

All along the proofs we will make use of the discrete Fourier
transform defined by 
\begin{equation*}
q_k=\frac{1}{\sqrt{2\pi}}\int_{-\pi}^\pi\hat q(\theta)\e^{\im
  k\theta}\d \theta\ ,
\end{equation*}
with 
\begin{equation*}
\hat q(\theta)=\frac{1}{\sqrt{2\pi}}\sum_{k\in\Z}q_k\e^{-\im
  k\theta}\ .
\end{equation*}
As usual the key property is that 
\begin{equation}
\label{delta1}
(\Delta q)^\wedge(\theta)=-(2-2\cos \theta)\hat
q(\theta)=-\left[4\sin^2\frac{\theta}{2}\right] \hat q(\theta)\ ,
\end{equation}
where
$$
(\Delta q)_k:=q_{k+1}+q_{k-1}-2q_k
$$
is the discrete Laplacian. 

\begin{lemma}
\label{dec.1} There exists $\epsilon_0$ s.t., if
$0<\epsilon<\epsilon_0$ then the operator $S^0_\epsilon(t)$ fulfills
\begin{equation}
\label{sti.3}
\norma{S^0_\epsilon(t)\xi}_{\bl^\infty}\sleq\frac{1}{\langle
  t\epsilon\rangle^{1/3}}\ .
\end{equation}\end{lemma}
\proof Rewrite \eqref{h.l.eq} as a second order equation in Fourier
coordinates, then it takes the form
\begin{equation}
\label{second}
\frac{d^2\hat q}{dt^2}(\theta)=-\nu(\theta)^2\hat q(\theta)\ ,\quad
\nu(\theta):= \sqrt{1+4{\epsilon}\sin^2\frac{\theta}{2}}\ ,
\end{equation}
whose solution is
$$
\hat q(\theta,t)=\hat q(\theta,0)\cos(\nu(\theta)t)+\frac{\hat
  p(\theta,0)}{\nu(\theta)} \sin (\nu(\theta)t)\ ,
$$
from which, returning to the space variables one gets
$$
q_k(t)=\frac{1}{\sqrt{2\pi}}\int_{-\pi}^\pi \left[\hat q(\theta,0)\cos(\nu(\theta)t)+\frac{\hat
  p(\theta,0)}{\nu(\theta)} \sin (\nu(\theta)t)\right]\e^{\im k\theta}
\d \theta 
$$
which is the linear combination of integrals of the form
\begin{equation}
\label{int.1}
\frac{1}{\sqrt{2\pi}}\int_{-\pi}^\pi\hat q(\theta,0)\e ^{\pm
  \im\nu(\theta)t+\im k\theta}\d \theta=
\sum_{j\in\Z}\frac{q_j(0)}{2\pi} \int_{-\pi}^\pi\e ^{\pm
  \im\nu(\theta)t+\im (k-j)\theta}\d \theta
\end{equation}
 and of the corresponding terms with $p$ instead of $q$. 

We now estimate the integral at r.h.s.~using the Van der Corput Lemma
(Lemma \ref{B.2} of the appendix). Writing $\varphi
(\theta,\rho):=\nu(\theta)+\rho\theta $, one has
$$
\sup_{k,j}\left|\int_{-\pi}^\pi  \e ^{\pm
  \im\nu(\theta)t+\im (k-j)\theta}\d \theta  \right|\leq
\sup_{\rho\in\R} \left|\int_{-\pi}^\pi  \e ^{
  \im \varphi(\theta,\rho)t}\d \theta\right|\ ,
$$ where, for definiteness, we choosed the sign +.  We split the
interval of integration $[-\pi,\pi]=I_1\cup I_2$ with
\begin{align*}
I_1:=\left[0,\frac{\pi}{8}\right]\cup
\left[\frac{3\pi}{8},\frac{5\pi}{8}\right]\cup
\left[\frac{7\pi}{8},\pi\right] 
\\
I_2:=\left[\frac{\pi}{8},\frac{3\pi}{8}\right]\cup
\left[\frac{5\pi}{8},\frac{7\pi}{8}\right] \ ,
\end{align*}
so that one has
\begin{align*}
\left|\varphi''(\theta,\rho)\right|=\left|\epsilon\cos2\theta+
O(\epsilon^2)\right|  \geq C\epsilon\ ,\quad \forall \rho\in\R\ ,\quad
\forall \theta\in I_1
\\
\left|\varphi'''(\theta,\rho)\right|=\left|-2\epsilon\sin2\theta+
O(\epsilon^2)\right|  \geq C\epsilon\ ,\quad \forall \rho\in\R\ ,\quad
\forall \theta\in I_2\ .
\end{align*}
Thus by Lemma \ref{B.2} one has 
$$
\sup_{\rho\in\R} \left|\int_{I_1}  \e ^{
  \im \varphi(\theta,\rho)}\d \theta\right|\sleq\frac{1}{|\epsilon
  t|^{1/2}}\ ,
$$
$$
\sup_{\rho\in\R} \left|\int_{I_2}  \e ^{
  \im \varphi(\theta,\rho)}\d \theta\right|\sleq\frac{1}{|\epsilon
  t|^{1/3}}\ ,
$$
from which, using also 
$$
|q_k(t)|\leq \norma{q(t)}_{\ell^2}\sleq \norma{\xi(0)}_{\bl^2}\sleq
\norma{\xi(0)}_{\bl^1}\ , 
$$ to control $t\to0$, one gets
$$
|q_k(t)|\sleq  \norma{\xi(0)}_{\bl^1}\min\left\{ 1,\frac{1}{|\epsilon
  t|^{1/2}}+\frac{1}{|\epsilon
  t|^{1/3}}  \right\}\sleq \frac{  \norma{\xi(0)}_{\bl^1}
  }{\langle\epsilon t\rangle^{1/3}}\ .
$$
Similarly, using
$$
\hat p(\theta,t)=\hat p(\theta,0)\cos(\nu(\theta)t)-{\hat
  p(\theta,0)}{\nu(\theta)} \sin (\nu(\theta)t)\ ,
$$
one gets the estimate of $|p_k(t)|$ and the proof of the Lemma.\qed

Next we need to establish weighted decay estimates. 

In the following we will denote by $B(\bl^2_s,\bl^2_{-s})$ the space
of bounded linear operators from $\bl^2_s$ to $\bl^2_{-s}$.

\begin{lemma}
\label{dec.2}
Let $s>5/2$ then one has
\begin{equation}
\label{s.pes}
\norma{S^0_\epsilon(t)\xi}_{\bl^2_{-s}}\sleq \frac{1}{\left\langle \epsilon
  t\right
\rangle^{3/2} } \norma{\xi}_{\bl^2_s}\ ,
\end{equation}
for all skew-symmetric sequence $\xi\in\bl^2_s$.   
\end{lemma}
\proof We follow closely the procedure of \cite{KKK} and use their
results (summarized in the appendix for the reader's convenience). 

First rewrite \eqref{h.l.eq} as $\im \dot\xi=A\xi$, where
\begin{equation*}
A:=\im\left[
\begin{matrix}
0 & -B
\\
\uno & 0
\end{matrix}\right]\ ,\quad B:= \uno-\epsilon \Delta\ ,
\end{equation*}
then it is easy to see that the spectrum $\sigma(A)$ of $A$ is given
by $\sigma(A)=I_+\cup I_-$ with $I_{\pm}:=\pm [1,\sqrt{1+4\epsilon}]$

An explicit computation shows that the resolvent
$R_A(\nu):=(A-\nu)^{-1}$ can be expressed in terms of the resolvent
$R_B$ of $B$ as follows
\begin{equation}
\label{Res.1}
R_A(\nu)=\left[
\begin{matrix}
\nu R_B(\nu^2) &- \im (\uno+\nu^2 R_B(\nu^2))
\\
\im R_B(\nu^2)&\nu R_B(\nu^2)
\end{matrix}
\right]\ ,\quad \nu\not \in \sigma(A) 
\end{equation}
 
Furthermore,remark that
\begin{equation}
\label{res.2}
R_B(\nu)= (1-\epsilon\Delta-\nu)^{-1}=
\frac{1}{\epsilon} \left(-\Delta-\frac{\nu-1}{\epsilon}\right)^{-1}=
\frac{1}{\epsilon}R_{-\Delta}\left(\frac{\nu-1}{\epsilon}\right)\ ,
\end{equation}
so that, from Lemma 3.1 of \cite{KKK} (see equation \eqref{Res.12}
below), the following limit exists in $B(\bl^2_s,\bl^{2}_{-s})$,
$s>1/2$
$$
R^{\pm}_A:=\lim_{\mu\to 0^+}R_A(\nu\pm\im \mu)\ , \nu\in
I_-\cup I_+\ .
$$
Let $\Gamma_{\pm}$ be closed curves enclosing $I_{\pm}$ respectively,
then by Cauchy theorem one has
$$
S^0_\epsilon(t)=\frac{1}{2\pi\im}\int_{\Gamma_+\cup\Gamma_-}\e^{-\im t
  \nu}R_A(\nu)\d \nu\ .
$$
We analyze the integral over $\Gamma_-$:
$$
\frac{1}{2\pi\im}\int_{\Gamma_-}\e^{-\im t
  \nu}R_A(\nu)\d \nu=\frac{1}{2\im\pi}\int_{-1}^{-\sqrt{1+4\epsilon}} \e^{-\im t
  \nu}\left[R_A^+(\nu)-R_A^-(\nu)\right]\d \nu\ ;
$$ making the change of variable $\nu=1+\epsilon\omega$ and exploiting
\eqref{res.2}, one gets that such quantity coincides with
\begin{align}
\label{int.3}
&\frac{\e^{-\im
    t}}{2\im\pi}\int_0^{-\frac{\sqrt{1+4\epsilon}-1}{\epsilon}}\e^{-\im\epsilon
t\omega} \times
\\ \nonumber
&
\left[
\begin{matrix}
(1+\epsilon\omega)[R^+_{-\Delta}(\varsigma(\omega))-R^-_{-\Delta}(\varsigma(\omega)
    )] &-\im(1+\omega\epsilon)^2 [R_{-\Delta}^+(\varsigma(\omega))-R_{-\Delta}^-(\varsigma(\omega)
    )]
\\
\im [R^+_{-\Delta}(\varsigma(\omega))-R^-_{-\Delta}(\varsigma(\omega)
    )] &(1+\epsilon\omega) [R^+_{-\Delta}(\varsigma(\omega))-R^-_{-\Delta}(\varsigma(\omega)
    )]
\end{matrix}
\right]\d \omega
\end{align}
where 
$$
\varsigma(\omega):=\frac{(1+\epsilon\omega)^2-1}{\epsilon}=
2\omega+\epsilon\omega^2 
\ .
$$
Using $\overline{R^+_{-\Delta}(\varsigma)}=R^-_{-\Delta}(\varsigma)$,
from which 
$$
\frac{R^+_{-\Delta}(\varsigma)-R^-_{-\Delta}(\varsigma
    )}{2\im}= \Im(R^+_{-\Delta}(\varsigma))\ ,
$$
one has that \eqref{int.3} coincides with
\begin{align}
\label{int.4}
\frac{\e^{-\im
    t}}{\pi}\int_0^{-\frac{\sqrt{1+4\epsilon}-1}{\epsilon}}\e^{-\im\epsilon
t\omega}
\left[
\begin{matrix}
1+\epsilon\omega & -\im(1+\epsilon\omega)
\\
\im & 1+\epsilon\omega
\end{matrix}
\right]
\\
\nonumber
\times
\left[
\begin{matrix}
\Im R^+_{-\Delta}(\varsigma(\omega)) & 0
\\
0& \Im R^+_{-\Delta}(\varsigma(\omega))
\end{matrix}
\right]\d \omega
\end{align}
Exploiting Lemma \ref{3.2KKK} one verifies that we are now in the
assumptions of Lemma \ref{JK}, which thus implies that the integral
\eqref{int.4} is bounded by a constant times $|\epsilon
t|^{-3/2}$. Treating in the same way the integral over $\Gamma_+$ one
gets the result. \qed

\begin{corollary}
\label{cor.dec}
Let $S_\epsilon(t)$ be the flow of the system with
Hamiltonian \eqref{d.h.l}, then, for any $s>5/2$, one has 
\begin{align}
\label{st.v.1}
\norma{S_\epsilon(t)\xi}_{\bl^2}\sleq \norma{\xi}_{\bl^2}\ '
\\
\label{st.v.2}
\norma{S_\epsilon(t)\xi}_{\bl^\infty}\sleq
\frac{\norma{\xi}_{\bl^1}}{\langle\epsilon t\rangle^{1/3}}
\\
\label{st.v.3}
\norma{S_\epsilon(t)\xi}_{\bl^2_{-s}}\sleq
\frac{\norma{\xi}_{\bl^2_s}}{\langle\epsilon t\rangle^{3/2}}\ .
 \end{align}
\end{corollary}

\subsection{Strichartz estimates}\label{sb.stri}

We first define the space-time norms which are needed in connections
with Strichartz inequalities.

The space $L^q_{\epsilon t}([0,T],\bl^r_s)$ is the space of the
functions $F:[0,T]\to \bl^r$ of class $L^q$ endowed by the norm 
\begin{equation}
\label{l.qepsi}
\norma{F}_{L^q_{\epsilon t}\bl^r_s} :=\left[\int_0^T\norma
    {F(t)}_{\bl^r}^q \epsilon \d t\right]^{1/q}=
\epsilon^{1/q}\norma{F}_{L^q_t\bl^r_s}\ ,
\end{equation}
where the last norm is defined in the usual way. In most cases we will
omit the indication of the interval of time and denote such a space
simply by $L^q_{\epsilon t}\bl^r$. 

We use the result of \cite{KT} to get the Strichartz estimates for our
model.

\begin{lemma}
\label{l.st.1}
Let $(q,r)$ and $(\tilde q,\tilde r)$ be admissible pairs, then the
flow $S_\epsilon(t)$ of \eqref{d.h.l} fulfills
\begin{align}
\label{St.1}
\norma{S_\epsilon(t)\xi}_{L^q_{\epsilon t}\bl^{r}}\sleq
\norma{\bl^2}
\\
\label{St.2}
\norma{\int_0^tS_\epsilon(t-\tau)F(\tau)\d \tau}_{L^q_{\epsilon
    t}\bl^r}\sleq \frac{1}{\epsilon}\norma {F}_{L^{\tilde q'}_{\epsilon
    t}\bl^{\tilde r'}} 
\end{align}
where $\tilde q' $ is such that $\frac{1}{\tilde q'}+\frac{1}{\tilde
  q}=1$ and similarly $\tilde r'$.
\end{lemma}

\proof Since $S_\epsilon(t)$ is not unitary with respect to the norm
of $\bl^2$ we first modify the norm suitably.  For $\xi\equiv
(p,q)\in \bl^2 $ we define
\begin{equation}
\label{n.b}
\norma{\xi}_{\bl^2_B}^2:= \langle p;p\rangle_{\ell^2}+\langle
q;Bq\rangle_{\ell^2}\ ,
\end{equation} 
where, as above, $B=\uno-\epsilon\Delta$ and, in the second scalar
product, one has to define $q_0\equiv 0$. It is immediate to verify
that in this metric $S^*(t)= S(-t)$. Furthermore the norm \eqref{n.b}
is equivalent to the standard norm of $\ell^2$. Moreover the norm
$\norma{ Bq}_{\ell^r_s}$ is equivalent to the norm
$\norma{q}_{\ell^r_s}$.

Before really starting with the proof remark that one has
$$
\norma{f\left(\frac{\cdot}{\epsilon}\right)}_{L^q_t}=\norma{f}_{L^q_{\epsilon
t}}\ ,
$$ and that $S\left(\frac{\cdot}{\epsilon}\right)$ is a group
fulfilling decay estimates independent of $\epsilon$. Thus Theorem 1.2
of \cite{KT} directly applies giving \eqref{St.1}. To get \eqref{St.2}
one has
\begin{align*}
\normaste{\int_0^tS_\epsilon(t-\tau)F(\tau)\d \tau }qr=
\normast{\int_0^{t/\epsilon}S_\epsilon\left(\frac{t}{\epsilon}-
  \tau\right)F(\tau)\d \tau }qr
\\
=\frac{1}{\epsilon}\normast{\int_0^{t/\epsilon}S_\epsilon\left(\frac{t-\tau'}{\epsilon}\right)F\left(\frac{\tau'}{\epsilon}\right)\d \tau' }qr
\\
\sleq \frac{1}{\epsilon}\normast{F\left(\frac{\cdot
  }{\epsilon}\right)}{\tilde q'}{\tilde
  r'}=\frac{1}{\epsilon}\normaste{F}{\tilde q'}{\tilde r'}\ .
\end{align*}
where the inequality is obtained by eq. (7) of \cite{KT}. \qed 

\begin{lemma}
\label{l.w.st}
Fix $s>5/2$ then, for any admissible pair $(q,r)$ one has
\begin{align}
\label{st.w.2}
&\norma{S_\epsilon(t)\xi}_{\bl^{\infty}_{-s}L^2_{\epsilon t}}\sleq
\norma{\xi}_{\bl^2}\ ,
\\
&\label{st.w.1}
\norma{\int_0^tS_\epsilon(t-\tau)F(\tau)\ d\tau}_{\bl^{\infty}_{-s}L^{2}_{\epsilon
t}} \sleq \frac{1}{\epsilon}\norma{F}_{\bl^{1}_sL^2_{\epsilon t}}\ ,
\\
&\label{st.w.3}
\norma{\int_0^tS_\epsilon(t-\tau)F(\tau)\ d\tau}_{\bl^{\infty}_{-s}L^{2}_{\epsilon
t}} \sleq \frac{1}{\epsilon}\norma{F}_{L^1_{\epsilon t}\bl ^2}\ ,
\\
&\label{st.w.4}
\normaste{\int_0^tS_\epsilon(t-\tau)F(\tau)\ d\tau}qr\sleq
\frac{1}{\epsilon}\norma{F}_{L^2_{\epsilon t}\bl ^2_s} 
\end{align}
\end{lemma}
\proof The proof is a minimal variation of Lemma 6 of \cite{KPS09}. We
begin by \eqref{st.w.1}. This is the equivalent of equation (27) of
\cite{KPS09}.  Since (27) is a consequence of the local decay
estimate, eq. \eqref{s.pes} with $\epsilon=1$ implies, by the
procedure of \cite{KPS09}, the validity of
\eqref{st.w.1} in the case $\epsilon=1$. The case with
$\epsilon\not=0$ is an immediate consequence of the same scaling
argument used in the proof of Lemma \ref{l.st.1}.

Equation \eqref{st.w.2} follows by the $TT^*$ argument when one
considers $T:\ell^2\to\bl^{\infty}_{-s}L^2_{\epsilon t}$. (Remark that
\eqref{st.w.2} is weaker then the corresponding equation in
\cite{KPS09}, namely (25).) 

Equations \eqref{st.w.3} and \eqref{st.w.4} follow from the previous
ones by repeating exactly the argument in the proof of Lemma 6 of
\cite{KPS09}. \qed

\subsection{Nonlinear estimates}\label{nonstri}

Here we prove the following Lemma:

\begin{lemma}
\label{non.lin}
Fix $\delta>1/2$, then there exists $\epsilon_\delta>0$ s.t., if
$0<\epsilon<\epsilon_\delta$ then the following holds true. Let
$(I(t),\alpha(t),\xi(t))$ be a solution of the Hamiltonian system
\eqref{H.tra} with initial datum $(I_0,\alpha_0,\xi_0)\in [\Delta_1-\frac{1}{2K_1},\Delta_2+\frac{1}{2K_1}]\times\toro\times\bl^2$ s.t.
\begin{align}
\label{non.rex}
\mu:=\norma{\xi_0}_{\bl^2}<\epsilon^\delta\ ,
\end{align}
then, for any admissible pair $(q,r)$ and any $s>5/2$ one has
\begin{align}
\label{rex.rex.1}
\norma{\xi}_{L^q_{\epsilon t}\bl^r}\sleq \mu\ ,
\\
\label{rex.rex.2}
\norma{\xi}_{\bl^\infty_{-s}L^2_{\epsilon t}}\sleq \mu \ .
\end{align}
Furthermore the limit $I_{\pm}:=\lim_{t\to\pm\infty}I(t)$ exists and
fulfills
\begin{equation}
\label{rex.rex.3}
|I_{\pm}-I_0|\sleq \frac{\mu^2}{\epsilon^{1/2}}\ .
\end{equation} 
\end{lemma}

\proof We proceed by ``induction'' as in \cite{GSNT04}: we are going to
prove that, if the solution fulfills
\begin{align}
\label{M.1}
\norma{\xi}_{L^q_{\epsilon t}([0,T],\bl^r)}\leq M_1 \mu\ ,
\\
\label{M.2}
\norma{\xi}_{\bl^\infty_{-s}L^2_{\epsilon t}[0,T]}\leq M_2 \mu\ ,
\end{align}
then it belongs to the interior of the domain of definition of 
the transformation of Thoeorem \ref{t.nor.for}, see item i), and
furthermore for a suitable choice of $M_1,M_2$ and for $\epsilon$
small enough, the inequalities \eqref{M.1} and \eqref{M.2} hold with
$M_1,M_2$ replaced by $M_1/2$ and $M_2/2$.

Using Duhamel formula rewrite the equation for $\xi$ in the form
\begin{align}
\label{var.xi}
\xi(t)=S_{\epsilon}(t)\xi_0+\int_0^t
S_\epsilon(t-\tau)X_{\xi}(I(\tau),\alpha(\tau), \xi(\tau))\d \tau
\\
\nonumber
+\int_0^t
S_\epsilon(t-\tau)\left[X_{\V}(\xi(\tau))\right]\xi\d \tau\ .
\end{align}
We begin by estimating the norm $L^q_{\epsilon t}\bl^r$ (where we
omitted the interval of time). The first term at r.h.s. is estimated
using \eqref{St.1} by
$$
\norma{S_{\epsilon}(t)\xi_0}_{L^q_{\epsilon t}\bl^r}\sleq
\norma{\xi_0}_{\bl^2}=\mu\ .
$$
For the second term, using \eqref{st.w.4} one has
\begin{align*}
\norma{\int_0^t
S_\epsilon(t-\tau)X_{\xi}(I(\tau),\alpha(\tau), \xi(\tau))\d
\tau}_{L^q_{\epsilon t}\bl^r}\sleq \frac{1}{\epsilon}
\norma{X_\xi(I,\alpha,\xi)}_{L^2_{\epsilon t}\bl^2_s}
\\
\sleq  \frac{1}{\epsilon}
\norma{X_\xi(I,\alpha,\xi)}_{L^2_{\epsilon t}\bl^\infty_{s'}} \sleq
\frac{1}{\epsilon} \epsilon^{3/2} \norma{\xi}_{L^2_{\epsilon
    t}\bl^{\infty}_{-s''}} \sleq
\epsilon^{1/2}\norma{\xi}_{L^2_{\epsilon t}\bl^2_{-s''}}\ ,
\end{align*}
where $s'>s+1/2$ and we used \eqref{sti.xi} for the third inequality,
which is valid for any $s''>0$. Using Lemma \ref{sp.temp} the last
quantity is smaller then 
$$
\epsilon^{1/2}\norma{\xi}_{\bl^{\infty}_{-s}L^2_{\epsilon t}}\leq
\epsilon^{1/2}M_2\mu \ ,
$$
provided $s''>s+1/2$.

For the third term one has, using \eqref{St.2},
\begin{align}
\nonumber
\norma{\int_0^t
S_\epsilon(t-\tau)\left[X_{\V}(\xi(\tau))\right]\xi\d
\tau}_{L^q_{\epsilon t}\bl^r}\sleq \frac{1}{\epsilon}
\norma{X_{\V}(\xi)}_{L^1_{\epsilon t }\bl^2}
\\
\label{sti.34}
\sleq \frac{1}{\epsilon} \norma{\xi}^7_{L^7_{\epsilon t}\bl^{14}}\sleq
\frac{\mu^{7}}{\epsilon}M_1^{7} \ , 
\end{align}
where we used, for $p=7$, the following inequalities 
\begin{align*}
\norma{X_{\V}(\xi)}_{L^1_{\epsilon t
  }\bl^2}=\int_0^T\left[\sum_{k\not=0}(V'(q_k(t)))^2
  \right]^{1/2}\epsilon\d t
\\
\sleq 
\int_0^T\left[\sum_{k\not=0}|q_k(t)|^{2p}
  \right]^{1/2}\epsilon\d t =\int_0^T\norma{q(t)}_{\bl^{2p}}^p
\epsilon\d t= \norma{\xi}_{L^p_{\epsilon t}\bl^{2p}}\ , 
\end{align*}
and the fact that $(p,2p)=(7,14)$ is an admissible pair.

Thus we have that the considered solution fulfills the inequality
\eqref{M.1} with $M_1$ replaced by $M_1/2$ if the following inequality
holds
\begin{equation}
\label{M.12}
1+M_2\epsilon^{1/2}+\frac{M_1^7\mu^{6}}{\epsilon}\leq \frac{M_1}{C}
\end{equation}
with a given large $C$. 

We estimate now $\norma{\xi}_{\bl^\infty_{-s}L^2_{\epsilon
    t}}$. Making again reference to equation \eqref{var.xi}, by
equation \eqref{st.w.2} one has
$$
\norma{S_\epsilon(t)\xi_0}_{\bl^\infty_{-s}L^2_{\epsilon
    t}}\sleq \norma{\xi_0}_{\bl^2}\sleq C_0\mu\ .
$$

Then, by \eqref{st.w.1} one has
\begin{align}
\nonumber
&\norma{\int_0^t
S_\epsilon(t-\tau)X_{\xi}(I(\tau),\alpha(\tau), \xi(\tau))\d
\tau}_{\bl^\infty_{-s}L^2_{\epsilon
    t}}
\\
\label{sti.33}
&\sleq
\frac{1}{\epsilon}\norma{X_\xi(I,\alpha,\xi)}_{\bl^{1}_sL^2_{\epsilon
    t}} \sleq \frac{1}{\epsilon}\norma{X_\xi(I,\alpha,\xi)}_{L^2_{\epsilon
    t} \bl^{2}_{s'}}
\end{align}
where $s'>s+1/2$ (the proof of the last inequality is almost identical
to the proof of Lemma \ref{sp.temp} and is
omitted). Eq. \eqref{sti.33} is estimated by using \eqref{sti.xi} and
lemma \ref{sp.temp}: 
$$
\eqref{sti.33}\sleq \frac{\epsilon^{3/2}}{\epsilon}\norma{\xi}
_{L^2_{\epsilon t}\bl^{2}_{-s''}}\sleq
\epsilon^{1/2}\norma{\xi}_{\bl^{\infty}_{-s}L^2_{\epsilon t}}\leq
M_2\epsilon^{1/2}\mu\ .
$$
The last term is estimated by eq.\eqref{st.w.3}, which gives, like in
\eqref{sti.34}, 
\begin{align*}
\norma{\int_0^t
S_\epsilon(t-\tau)\left[X_{\V}(\xi(\tau))\right]\xi\d
\tau}_{\bl^\infty_{-s}L^2_{\epsilon
    t}  }\sleq \frac{1}{\epsilon}\norma{X_{\V}(\xi)}_{L^1_{\epsilon
    t}\bl ^2}
\\
\sleq  \frac{\mu^7}{\epsilon}M_1^7\ ,
\end{align*}
so that the considered solution fulfills the inequality
\eqref{M.2} with $M_2$ replaced by $M_2/2$ if the following inequality
holds
\begin{equation}
\label{M.22}
1+M_2\epsilon^{1/2}+\frac{M_1^7\mu^{6}}{\epsilon}\leq \frac{M_2}{C}
\end{equation}
with a given large $C$.
 
Now it is clear that both \eqref{M.12} and \eqref{M.22} are fulfilled if
$M_1$ and $M_2$ are chosen strictly larger then $C$, with $C$ the
constant in \eqref{M.12} and \eqref{M.22}, and $\epsilon$ is small
enough. In particular this implies that also $\mu^6/\epsilon$ is
small.

Remark that from these inequalities it also follows that
$\norma{\xi(t)}_{\bl^2}<\sqrt \epsilon/2K_1$, so that $\xi$ is in the
domain of validity of the normal form.

Concerning $I$, one has
$$
I(t)=I_0+\int_0^tX_I(\zeta(\tau))\d \tau\ .
$$
One has 
\begin{align}
\label{M.32}
\int_0^t\left|X_I(\zeta(\tau))\right|\d \tau\sleq
\frac{\epsilon^{1/2}}{\epsilon}
\int_0^t\norma{\xi(\tau)}^{2}_{\bl^2_{-s}}\epsilon \d
\tau=\frac{1}{\epsilon^{1/2}} \norma{\xi}^2_{L^2_{\epsilon
    t}\bl^2_{-s'}}
\\
\label{M.31}
\sleq \frac{1}{\epsilon^{1/2}}\norma{\xi}^2_{\bl^\infty_{-s}L^2_{\epsilon
    t}}\sleq \frac{\mu^2}{\epsilon}M^2_2\ ,
\end{align}
which implies $|I(t)-I(0)|\sleq \epsilon^{2\delta-1}$ and therefore
also $I(t)$ is close to $I_0$ and therefore, if $\epsilon$ is small
enough $I(t)\in[\Delta_1-\frac{3}{4K_1},\Delta_2+\frac{3}{4K_1})$
  which is contained in the domain of validity of the normal form.

This allows to extend the estimates to $T=\infty$.  From \eqref{M.31}
follows that the integral \eqref{M.32} converges, and therefore the
limit of $I(t)$ exists and \eqref{rex.rex.3} holds.\qed

\noindent {\it{Proof of Theorem \ref{breather}}.} First, the existence
of the breather and item i) are a consequence of theorem
\ref{t.nor.for} (see Remark \ref{cor.bre}).

Item ii.1 follows immediately by defining $\I$ as the action variable
$I$ in the coordinates introduced by the normal form theorem. 

Finally, to get \eqref{asybreather}, remark that, in the coordinates
introduced by Theorem \ref{t.nor.for} 
$$
d_{\bl^{r}}(\gamma_\epsilon(\I(t));\zeta(t))=\norma{\xi(t)}_{\bl^r}\ ,
$$ then \ref{asybreather} follows from \eqref{rex.rex.1} and the fact
that the canonical transformation $T$ is Lipschitz in the $\bl^r$
metric (see Remark \ref{lr.analit}), and therefore only multiplies
distances by a number (which is of order 1 in our case).\qed

\appendix
\section{Technical lemmas for the normal form}\label{technical}

We begin by the different estimates involved in Lemma
\ref{sti.varie}. 

The estimate \eqref{4.3.1} of the Poisson brackets coincides with
that given in \cite{BG93}, lemma 5.2. For the sake of completeness we
repeat here the argument of that paper.
\begin{lemma}
\label{4.1}{Let $g\in\A_d$ and $f\in\Sc_{d}$ be two functions
  with analytic vector field; then for any $d_1<1-d$, $ \poisson
  g{f}\in\Sc_{d+d_1}$ satisfies the inequality \eqref{4.3.1}.  }
\end{lemma}

\proof First remark that 
\begin{equation}
\label{lie}
X_{\poisson f{g}}=[X_f;X_{g}]= \di X_f\, X_{g}- \di X_{g}\,
X_{f}\ .
\end{equation}
Using Cauchy estimate one immediately has
that, on $\G^{-}_{d+d_1}$, the norm of ${\di X_f}$ as a
linear operator from $\Ph^-$ to $\Ph^+$ is smaller then $\normas d
f/d_1$. It follows that the norm $\normas{d+d_1}.$ of
first term of \eqref{lie} is bounded by
$$
\frac1{d_1}\normas
d f \normad d{g}\ .
$$
The second term is bounded in a similar way getting the thesis.\qed
\vskip20pt

\begin{lemma}
\label{parti}
Let $f\in\Sc_d$ be a function with analytic vector field; let
$0<d_1<1-d$, then the estimates \eqref{sti.parti.21} and
\eqref{stime.parti.111} hold.
\end{lemma}
\proof The estimate is trivial for the $\null^{(0)}$ component. Indeed
the vector field of $f^{(0)}$ coincides with the value at $\xi=0$ of
the components $(I,\alpha)$ of the vector field of $f$.

We come to the estimate of the vector field of
$f^{(1)}(I,\alpha,\xi)\equiv \di_\xi f(I,\alpha,0)\xi$. Remark that one
has
\begin{align}
\label{sti.parti.2}
\left[X_{f^{(1)}}\right]_{\xi}(I,\alpha)=\left[X_{f}\right]_{\xi}(I,\alpha,0)
\ ,
\\
\label{sti.parti.3}
\left[X_{f^{(1)}}\right]_{x}(I,\alpha,\xi)=\di_{\xi}
\left[X_{f}\right]_{x}(I,\alpha,0) \xi
\ .
\end{align}
The estimate of \eqref{sti.parti.2} is straightforward. Concerning the
estimate of \eqref{sti.parti.3}, remark that, by Cauchy inequality one
has
\begin{equation}
\label{sti.p.4}
\norma{\di_{\xi}
\left[X_{f}\right]_{x}(I,\alpha,0)}\leq
\frac1{1-d}\normas{d}{f}\ ,
\end{equation}
and that, on $\G^-_d$ the norm \eqref{3.13} of $\xi$ is smaller
then $R(1-d)$. Thus one gets also the second of \eqref{sti.parti.21}.

We come to the estimate of $f^{(2)}$. The components of its vector field
 are remainders of Taylor expansions truncated at suitable order of
the components of
$X_f$. In particular the term of higher order is in the $x$
components. From standard formulae of the remainder of Taylor
expansions one has
$$
\left[
  X_{f^{(2)}}\right]_x(I,\alpha,\xi)=\int_0^1(1-s)\di^2_\xi\left[X_f\right]_x
(I,\alpha, s\xi)(\xi,\xi)\di s\ ;
$$
using Cauchy estimate to estimate the norm of the second differential
one gets that the argument of the integral, in
$\G^-_{d+d_1}$ is estimated by
$$
\frac{2}{R^2d_1^2}R\normas d f [R(1-d-d_1)]^2\ , 
$$ 
which, integrating and dividing by $R$ in order to get the norm
$\normas{d+d_1}.$ gives the result.
\qed

\begin{lemma}
\label{poi.part}
Let $f=f^{(1)}\in\Sc_ d$ and $g=g^{(2)}\in \A_ d$. Then
\eqref{poi.part.1.1} holds.
\end{lemma}
\proof First remark that, denoting by
$g_2(I,\alpha,\xi):=[\di^2_\xi g(I,\alpha,0)](\xi,\xi)$ the part of
$g^{(2)}$ homogeneous of degree 2, one has
$$
\poisson{f^{(1)}}{g^{(2)}}^{(1)}=\poisson{f^{(1)}}{g_2}^{(1)}\ .
$$ So we first study $g_2$. Remark that, by a procedure similar to the
one used in the proof of lemma \ref{parti}, one has
\begin{equation}
\label{z2}
\normad d{g_2}\leq \normad d{g^{(2)}}\ .
\end{equation}

Denote $ B(I,\alpha):=J^{-1}\di_\xi
\left[X_{g^{(2)}}\right]_\xi(I,\alpha,0)$, where $J$ is
the Poisson tensor, then one has 
\begin{equation}
\label{poi.p.2}
g_2(I,\alpha,\xi)=\frac{1}{2} \left\langle \xi;B(I,\alpha)\xi
\right\rangle \ ,\quad \norma{B}\leq
\frac{1}{1- d}\normad d{g^{(2)}}\ ,
\end{equation}
and furthermore $B$ is symmetric. The considered norm of $B$ is the
maximum between the norm as an operator from $\Ph^+$ to itself and as
an operator from $\Ph^-$ to itself.

So one has
\begin{align}
\label{poi.p.4}
\poisson{f^{(1)}}{g^{(2)}}^{(1)}
= \left\langle
f^1(I,\alpha),JB(I,\alpha)\xi   \right\rangle=-\left\langle
B(I,\alpha)Jf^1(I,\alpha),\xi   \right\rangle\ . 
\end{align}
From this formula one has that the $\xi$ component of the vector
field, given by $-JB(I,\alpha)f^1(I,\alpha)$ is actually estimated
by \eqref{poi.part.1.1}. 

We come to the $\alpha$ component of the vector field: it is given by 
$$
-\left\langle
\der B IJ f^1,\xi   \right\rangle-\left\langle
\der {f^1}I,JB\xi   \right\rangle\ .
$$ We start by estimating the second term. To this end remark that
$$\der {f^1}I(I,\alpha)=
\nabla_\xi\left[X_{f^{(1)}}\right]_\alpha(I,\alpha,0)\ ,$$
so that, exploiting the fact that $\bl^+$ is the dual of $\bl^-$, 
its norm \eqref{3.13} can be
bounded using Cauchy inequality:
$$
\snorma{\der{f^1}I(I,\alpha)}_+\leq
\frac{1}{R(1-d)}\sup_{\G^-_d}\left|\left[X_{f^{(1)}}\right]_\alpha\right|
\leq \frac{1}{1-d}\normas d{f^{(1)}}R_\alpha\ .
$$
Using also the
estimate \eqref{poi.p.2} one thus gets 
\begin{align*}
\frac{1}{R_\alpha}\left|\left\langle
\der {f^1}I,JB\xi   \right\rangle\right|\leq \frac1
     {(1- d)}\normas d{f^{(1)}}\norma
     B R(1- d)
\\
\leq
R\normas d
     {f^{(1)}} \frac{1}{1- d}\normad  d{g^{(2)}}\ ,
\end{align*}
dividing by $R$ one gets the wanted estimate.  We now estimate the
term involving the derivative of $B$. Remark first that one has
$$
\left[X_{g_2}\right]_\alpha=\left\langle \xi;\der BI\xi\right\rangle\ ,
$$
so that the norm of $\der BI$ as an operator from $\Ph^-$ to $\Ph^+$
is estimated by 
$$
\norma{\der BI}\leq
\frac2{R^2(1- d)^2}\sup_{\G^-_ d}\norma{\left[X_{g_2}\right]_\alpha}\leq
\frac{2}{R(1- d)^2}\normad d{g_2}\ ,
$$
from which, on $\G_d^-$,
\begin{align*}
\frac{1}{R_\alpha}\left|\left\langle
\der B IJ f^1,\xi   \right\rangle  \right|\leq
\frac{2}{R(1- d)^2}\normad d{g_2}
R(1- d)\sup_{\G_d^-}\norma{\left[X_{f^{(1)}} \right]_\xi}_+
\\
\leq \frac{2R}{(1- d)}\normad d{g_2}\normas d{f^{(1)}}
\ .
\end{align*}
Collecting the results the thesis follows.\qed

\begin{lemma}
\label{tra} Let
 $\chi\in\Sc_d$, with $0\leq d_1<1-d$
 and let $f\in\A_ d$; fix $0<
d_1<(1- d)$, assume $\normas d {\chi}\leq d_1/3$, then, for
$|t|\leq1$, one has
$$
\normad{ d+ d_1}{f\circ \Phi_\chi^t}\leq 
\left(1+\frac{3}{ d_1}\normas d{\chi}\right)
\normad d{f}\ .
$$
If $f\in\Sc_ d$ then the same estimate holds in the norm $\normas..$.
\end{lemma}
\proof In this proof we omit the index $\chi$ from $\Phi$. First
remark that, since $\Phi^t$ is a canonical transformation one has
\begin{equation}
\label{c.T}
 X_{f\circ \Phi^t}(\zeta)=
\di\Phi^{-t}(\Phi^t(\zeta))X_{f(\Phi^t(\zeta))}\ ,
\end{equation}
from which
$$
X_{f\circ \Phi^t}(\zeta)=\left(\di \Phi^{-t}(\Phi^t(\zeta))-\uno
\right)
X_{f(\Phi ^t(z))} +X_{f(\Phi^t(z))}\ .
$$
We first estimate $X_{f\circ\Phi^t}$ in $\Ph^+$.
To estimate the first term fix $\bar d:= d_1/3$; we have
\begin{align}
\label{c.T.1}
\sup_{\zeta\in\G^+_{3\bar  d}}\norma{\di \Phi^{-t}(\Phi^t(\zeta))-\uno }\leq 
\sup_{\zeta\in\G^+_{2\bar  d}}\norma{\di \Phi^{-t}(\zeta)-\uno} 
\\
\nonumber
\leq \frac1{\bar  d}
\sup_{\zeta\in\G^+_{\bar  d}}\norma{\Phi^{-t}(\zeta)-\zeta}_+\leq 
\frac1{\bar  d}\normas d{\chi}\ ,
\end{align}
where the differential of $\Phi^{-t}(\zeta)$ is considered as an
operator from $\Ph^+$ to $\Ph^+$.  Going back to $ d_1$, adding the
trivial estimate of the second term and the estimate in $\Ph_-$, one
gets the thesis.\qed

\noindent{\it Proof of Lemma \ref{l.2}}. First define $\bar
d:=\frac{d}{2(N+1)}$. Using \eqref{4.3.1} $l$-times, one has 
$$
\normas{d+l\bar d}{f_{(l)}}\leq \left(\frac{2}{\bar d}\normas{d}\chi \right)^l
\normad d f\ . 
$$
By Lemma \ref{tra} one has
$$
\normas{d+d_1}{f_{(N+1)}\circ\Phi^t}\leq \left(1+\frac{6}{d_1}\normas
d\chi \right)\left(\frac{4(N+1)}{d_1}\normas d\chi\right)^{N+1}
\normad d f\ ,
$$
which gives the
thesis. 
\qed

\noindent
{\it Proof of lemma \ref{sol.homo}}. We start by $\chi^{(0)}$. It is well known (see
e.g. \cite{BG93}) that defining 
\begin{equation}
\label{chi0}
\chi^{(0)}(I,\alpha):=\frac{1}{2\pi\omega(I)}
\int_0^{2\pi}t\Psi^{(0)}(I,\alpha+t)\di t\ ,
\end{equation}
it solves the equation 
$$
\poisson{H_{lin}}{\chi^{(0)}}= \Psi^{(0)}\ .
$$
Then one has
\begin{equation}
\label{xchi0}
X_{\chi^{(0)}}(I,\alpha):=\frac{1}{2\pi\omega(I)}
\int_0^{2\pi}tX_{\Psi^{(0)}}(I,\alpha+t)\di t\ ,
\end{equation}
from which the estimate \eqref{sti.homo} immediately follows. 

We now study the equation
\begin{equation}
\label{homo.11}
\poisson{H_{lin}}{\chi^{(1)}}= \Psi^{(1)}\ ;
\end{equation}
inserting the decomposition \eqref{p.1}, i.e. writing
\begin{equation}
\label{psi.1}
\Psi^{(1)}(I,\alpha,z,w)=\langle \Psi_{z}^1;z \rangle+\langle
\Psi^1_{w};w \rangle \ ,
\end{equation}
and similarly for $\chi$, one gets that \eqref{homo.11} is equivalent
to the couple of equations
\begin{align}
\label{homo.2}
\omega(I)\der{\chi^1_z}\alpha-\im \chi^1_z=\Psi^1_z\ ,\quad 
\omega(I)\der{\chi^1_w}\alpha+\im \chi^1_w=\Psi^1_w\ .
\end{align}
Let's focus on the second one. Component wise this is an ordinary
differential equation in the independent variable $\alpha$, which can
be easily solved by Duhamel formula. Imposing the solution to be
periodic of period $2\pi$ in $\alpha$ one gets a unique solution given
by
\begin{equation}
\label{sol.homo.3}
\chi^1_w(I,\alpha)=\frac{1}{\omega\left(\e^{\im2\pi\frac{1}{\omega}}-1
\right)}\int_0^{2\pi}\e^{\im2\pi\frac{1}{\omega}s}\Psi^1_w(I,\alpha+s)\di
s \ .
\end{equation}
For $\chi^1_z$ one gets an identical formula with $-1$ in place of
$1$. From \eqref{sol.homo.3} one gets identical formulae for the
functions $\langle\chi^1_w,w\rangle$, $\langle\chi^1_z,z\rangle$ and
for their Hamiltonian vector fields. Inserting the corresponding
estimates of the vector field of $\langle\Psi^1_w,w\rangle$ and of the
other functions and computing the integrals one gets the thesis. \qed

\section{Technical lemmas for the dispersive estimates}\label{techdis}

\begin{lemma}
\label{B.2}
Let $\varphi(\theta)$ be a function of class $C^k(a,b)$, $k\geq 2$
and assume that
$$
\left|\varphi^{(k)}(\theta)\right|\geq \delta_k>0\ ,\quad \forall
\theta\in (a,b)\ ,
$$
then there exists $c_k$ s.t.
\begin{equation}
\label{vander}
\left|\int_a^b\e^{\im\lambda\varphi (\theta)}\d
\theta\right|\leq\frac{c_k}{|\lambda\delta_k|^{1/k}}\ .
\end{equation}
\end{lemma}
The proof is a minor variant of the proof of Proposition 2 p.332 of
\cite{Stein} (Van der Corput Lemma), and is omitted.

We now recall the properties of $-\Delta$ and in particular the
Puiseaux expansion for $R_{-\Delta}$ proved in \cite{KKK} and
specialize it to skew symmetric sequences. 

By using an explicit computation and the Cauchy formula for the
computation of integrals \cite{KKK}, proved the following lemma:

\begin{lemma}
\label{KKK} (2.1 of \cite{KKK})For $\tilde \nu\in\C-[0,4]$ the Kernel
of the resolvent of $-\Delta$ is given by
\begin{equation}
\label{Res.11}
R_{-\Delta}(\tilde \nu,j,k)=-\im \frac{\e^{\im
    \theta(\tilde\nu)|j-k|}}{2\sin(\theta(\tilde\nu))}
\end{equation}
where $\theta(\tilde\nu) $ is the unique solution of the equation
\begin{equation*}
2-2\cos \theta=\tilde \nu
\end{equation*}
in the domain $\{-\pi\leq \Re\theta\leq \pi\ ;\ \Im\theta<0   \}$.
\end{lemma}
By this we mean that 
\begin{equation*}
(R_{-\Delta}(\tilde \nu)q)_j=\sum_{k}R_{-\Delta}(\tilde \nu,j,k)q_k\ .
\end{equation*}
\begin{corollary}
\label{KKK.1} Equation \eqref{res.2} holds for $\nu\in\C-[1,1+4\epsilon]$.
\end{corollary}

In Lemma 3.1 of \cite{KKK}, by direct computation of the limit of
\eqref{Res.11}, it is shown that the limit
\begin{equation}
\label{Res.12}
\lim_{\epsilon\to 0^+}R_{-\Delta}(\tilde\nu \pm\im
\epsilon)=R^{\pm}_{-\Delta}(\tilde \nu)\ ,\quad \tilde \nu\in (0,4)
\end{equation}
exists in $B(\ell^2_s,\ell^{2}_{-s})$ for all $s>1/2$, and this
implies a similar result for $R_B$.

Remark that the term proportional to $|\tilde \nu|^{-1/2}$ in our case
is missing. This is due to the fact that its coefficient is
proportional to $\sum_lq_l$, which vanishes for skewsymmetric
sequences.

Then we need the following lemma, which is a particular case of Lemma
3.2 of \cite{KKK}.

\begin{lemma}
\label{3.2KKK}
Let $s>3/2$, then for $\tilde \nu\to 0$ one has the following
asymptotic expansion, valid for {\bf skew-symmetric sequences} $q$
\begin{equation}
\label{poise.1}
\left[R^{\pm}_{-\Delta_1}(\tilde
  \nu)q\right]_k=-\frac{1}{2}\sum_{l}\left|k-l\right| q_l+r(\tilde
\nu)q\ ,
\end{equation}
where $\norma{r(\tilde
  \nu)}_{B(\ell^2_s,\ell^2_{-s})}=O(|\tilde\nu|^{1/2})$ and, for
$s>\frac{1}{2}+i$, $i\geq 1$ one has
\begin{equation}
\label{poise.2}
\norma{\frac{d^{i}}{d\tilde\nu^i}R^{\pm}_{-\Delta}}_{B(\ell^2_s,\ell^2_{-s})}
=O(|\tilde\nu|^{\frac{1}{2}-i})\ .
\end{equation}  
A similar expansion holds for $\tilde \nu\to 4$.
\end{lemma}

We also need the following Lemma by Jensen-Kato

\begin{lemma}
\label{JK}
Let $\Bc$ be a Banach space, and let $F\in C^2((0,a),\Bc)$, assume
$$
F(0)=F(a)=0\ ,\quad \norma{\frac{d^{i}}{d\tilde\nu^i}F(\tilde
  \nu)}_{\Bc} \leq C |\tilde \nu|^{\frac{1}{2}-i}\ ,\quad \tilde
\nu\to 0\ ,\quad i=1,2
$$ 
then for any $|t|>1$ one has 
\begin{equation}
\label{JK.1}
\left|\int_0^a\e^{\im t\tilde\nu}F(\tilde \nu)\d \tilde \nu\right|\leq
\frac{C}{|t|^{3/2}} 
\end{equation}
\end{lemma}
For the proof see \cite{JK} (see also \cite{PS08}).

\begin{lemma}
\label{sp.temp}
One has
\begin{align*}
\norma q_{L^2_{t}\ell^\infty_{-s}}\sleq \norma
q_{\ell^\infty_{-s'}L^2_t}\ ,\quad \forall s>s'+\frac{1}{2}\ .
\end{align*}
\end{lemma}
\proof One has
\begin{align*}
\norma{q}^2_{L^2_t\ell^\infty_{-s}}=\int\left[\sup_n\langle
  n\rangle^{-s}|q_n(t)| \right]^2\d t
\\
\leq \int\sum_n\langle n\rangle^{-2s}|q_n(t)|^2\d t= \sum_n \langle
n\rangle^{-2s} \int \left| q_n(t)\right|^2 \d t
\\
\leq \left[\sum_{n}\langle n\rangle^{-2(s-s')}  \right]\sup_n
\left[\left\langle n\right\rangle^{-2s'}\int \left|q_n(t)\right|^2\d t\right]
\\
=C\norma q^2_{\ell^\infty_{s'}L^2_t}
\end{align*}
\qed


\bibliographystyle{amsalpha}

\providecommand{\bysame}{\leavevmode\hbox to3em{\hrulefill}\thinspace}
\providecommand{\MR}{\relax\ifhmode\unskip\space\fi MR }
\providecommand{\MRhref}[2]{%
  \href{http://www.ams.org/mathscinet-getitem?mr=#1}{#2}
}
\providecommand{\href}[2]{#2}

\end{document}